\documentclass[english,twoside,11pt]{article}

\usepackage[german,english]{babel}
\usepackage{amsmath}
\usepackage{amssymb,latexsym}
\usepackage{moreverb}
\usepackage{url}
\usepackage{graphicx}
\usepackage{longtable}
\usepackage{lscape}
\usepackage{booktabs}
\usepackage{tabularx}
\usepackage{psfrag}
\usepackage{multicol}
\usepackage{paralist}

\newtheorem{deff}{Definition}

\newtheorem{prop}[deff]{Proposition}
\newtheorem{thm}[deff]{Theorem}

\newtheorem{cor}[deff]{Corollary}
\newtheorem{conj}[deff]{Conjecture}

\makeatletter

\title{Equivelar and $d$-Covered Triangulations of Surfaces. I}

\author{\Large Frank H.~Lutz, Thom Sulanke, Anand K.~Tiwari, Ashish K.~Upadhyay}

\date{}

\begin{document}

\selectlanguage{english}

\maketitle

\vspace{2mm}

\begin{abstract}
We survey basic properties and bounds for $q$-equivelar and $d$-covered triangulations
of closed surfaces. Included in the survey is a list of the known sources for $q$-equivelar 
and $d$-covered triangulations. We identify all orientable and non-orientable surfaces $M$
of Euler characteristic $0>\chi(M)\geq -230$ which admit non-neighborly $q$-equivelar
triangulations with equality in the upper bound
$q\leq\Bigl\lfloor\tfrac{1}{2}(5+\sqrt{49-24\chi (M)})\Bigl\rfloor$.
These examples give rise to $d$-covered triangulations
with equality in the upper bound $d\leq2\Bigl\lfloor\tfrac{1}{2}(5+\sqrt{49-24\chi (M)})\Bigl\rfloor$.
A~generalization of Ringel's cyclic $7{\rm mod}12$ series of neighborly orientable triangulations \cite{Ringel1961}
to a two-parameter family of cyclic orientable triangulations $R_{k,n}$, $k\geq 0$, $n\geq 7+12k$, 
is the main result of this paper. In particular, the two infinite subseries $R_{k,7+12k+1}$ and $R_{k,7+12k+2}$, $k\geq 1$,
provide non-neighborly examples with equality for the upper bound for $q$ as well as derived examples 
with equality for the upper bound for $d$.

\end{abstract}

\section{Introduction and Survey}

A triangulation (as a finite simplicial complex) $K$ of a (closed) surface $M$ is \emph{equivelar} 
if all vertices of the triangulation have the same vertex-degree $q$. 
The \emph{$f$-vector} (or \emph{face vector}) of~$K$ is the vector $f(K)=(f_0(K),f_1(K),f_2(K))$, 
where $f_0(K)$, $f_1(K)$, and $f_2(K)$ are denoting the numbers of vertices, edges, 
and triangles of $K$, respectively. For brevity, we write $f=(n,f_1,f_2)$ 
for the $f$-vector of a triangulation, with $n=f_0$ denoting the numbers of vertices.
If the vertices of an equivelar triangulation $K$ of a surface $M$ have degree $q$,
then the triangulation is called \emph{$q$-equivelar} or \emph{equivelar of type $\{3,q\}$}.
For $q$-equivelar triangulations it follows by double counting of incidences between vertices and edges 
as well as between edges and triangles that
\begin{equation}
nq=2f_1=3f_2,
\end{equation}
which yields $f=(n,\frac{nq}{2},\frac{nq}{3})$. Moreover, by Euler's equation, 
\begin{equation}\label{eq:equitri}
\chi(M)=n-f_1+f_2=n-\frac{nq}{2}+\frac{nq}{3}=\frac{n(6-q)}{6},
\end{equation}
or equivalently,
\begin{equation}\label{eq:q}
q=6-\frac{6\chi(M)}{n}.
\end{equation}
Since $q$ is a positive integer, $n$ has to be a divisor of $6|\chi(M)|$ if $\chi(M)\neq 0$.
In particular, a~surface~$M$ of Euler characteristic $\chi(M)\neq 0$ has only finitely many equivelar triangulations;
see for example \cite{SulankeLutz2009} for a list of the possible values of $(n,q)$ for
equivelar triangulations of surfaces $M$ with $\chi(M)\geq -10$.

\begin{figure}
\begin{center}
\includegraphics[width=.14\linewidth]{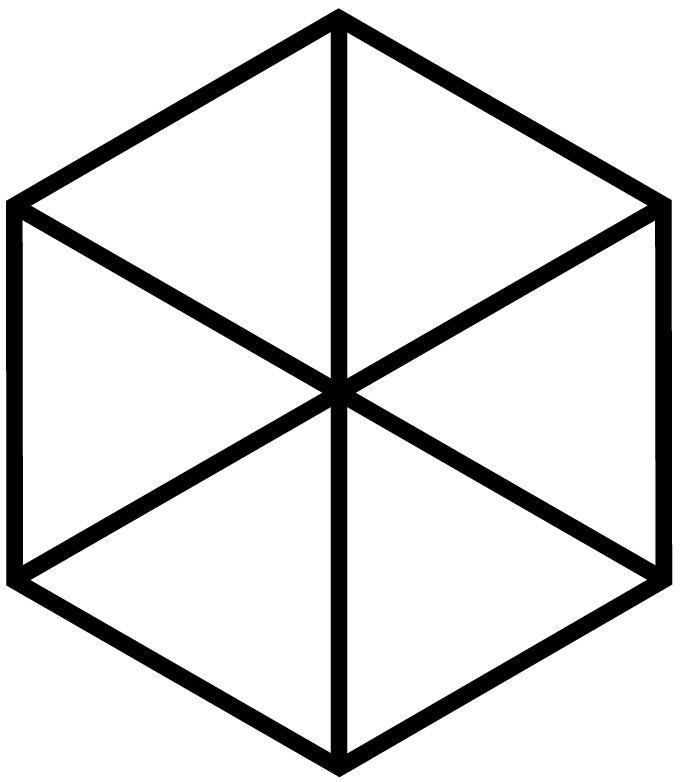}\hspace{20mm}\includegraphics[width=.14\linewidth]{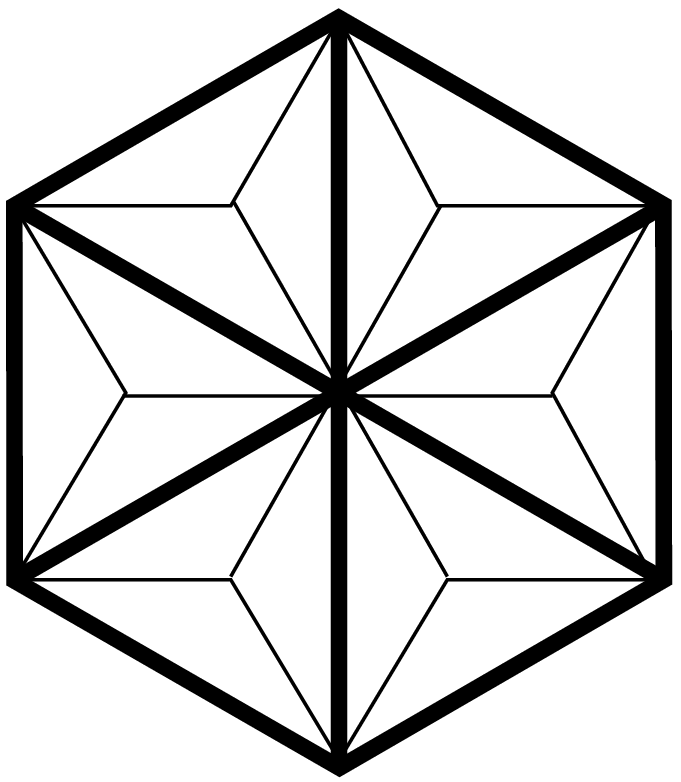}
\end{center}
\caption{A vertex-star of degree six and its subdivision.}
\label{fig:subdivision}
\end{figure}

A triangulation $K$ of a surface $M$ is \emph{$d$-covered}, if at least one vertex of each edge of $K$ has degree~$d$.
The study of $d$-covered triangulations was initiated by Negami and Nakamoto in \cite{NegamiNakamoto2001}. 
In particular, a $d$-covered triangulation $K$ is covered by triangulated disks, which are vertex-stars of vertices 
of degree $d$. Any $q$-equivelar triangulation clearly is a $q$-covered triangulation. If we \emph{stack} all the triangles 
of a $q$-equivelar triangulation, i.e., if we place for each triangle a new vertex inside the triangle and connect 
the new vertex with the three vertices of the triangle by an edge each, then the resulting triangulation is $2q$-covered \cite{NegamiNakamoto2001}; 
see Figure~\ref{fig:subdivision} for the respective subdivision of a vertex-star of degree~six. 

For an arbitrary triangulation $K$ of a surface $M$ of Euler characteristic $\chi(M)$, the $f$-vector of~$K$ can
(by double counting, $2f_1=3f_2$, and by Euler's equation, $\chi(M)=n-f_1+f_2$)
be expressed as $f=(n,3(n-\chi(M)),2(n-\chi(M)))$. Since $K$ has at most $f_1\leq\binom{n}{2}$ edges,
it follows that 
\begin{equation}\label{eq:Heawood}
n\geq\Bigl\lceil\tfrac{1}{2}(7+\sqrt{49-24\chi (M)})\Bigl\rceil,
\end{equation}
which is Heawood's bound~\cite{Heawood1890} from 1890. A triangulation $K$ is called \emph{neighborly}
if $f_1=\binom{n}{2}$, i.e., if $K$ has a complete $1$-skeleton. Neighborly triangulations are equivelar
with $q=n-1$. Moreover, we have equality $3(n-\chi(M))=f_1=\binom{n}{2}$ for neighborly triangulations,
or, equivalently
\begin{equation}\label{eq:chi_equiv}
\chi(M)=\frac{n(7-n)}{6},
\end{equation}
which implies $n\equiv 0,1,3,4\,{\rm mod}\,6$, or equivalently, $n\equiv 0,1\,{\rm mod}\,3$, where $n\geq 4$.

\begin{thm} \label{thm:Jungerman_Ringel}
{\rm (Jungerman and Ringel \cite{JungermanRingel1980}, \cite{Ringel1955})}
Let $M$ be a surface different from the orientable surface of genus $2$, 
the Klein bottle, and the non-orient\-able surface of genus $3$,
then there is a triangulation of $M$ on\, $n$ vertices if and only if
\begin{equation*}
n\geq\Bigl\lceil\tfrac{1}{2}(7+\sqrt{49-24\chi (M)})\Bigl\rceil,
\end{equation*}
with equality if and only if the triangulation is neighborly.
For the three exceptional cases, one extra vertex has to be added to the lower bound; cf.~{\rm \cite{Huneke1978}}.
\end{thm}

Let $G=\mbox{\rm Sk}_1(K)$ be the $1$-skeleton of $K$, $\mbox{\rm deg}_{\,G}(v)$ the degree of a vertex $v$ in the graph $G$,
and $\delta(G)=\displaystyle\mbox{\rm min}_{\,v\in G\,}\mbox{\rm deg}_{\,G}(v)$ the minimum degree in $G$. 
For any triangulated surface $K$ with $n$ vertices it trivially holds that $3\leq \delta(G)\leq n-1$.
Moreover, $n\delta(G)\leq\sum_{\,v\in G\,}\mbox{\rm deg}_{\,G}(v)=2f_1$, 
with the equality by double counting of vertex-edge incidences.
By combining $\delta(G)\leq \frac{2f_1}{n}$ with $f_1=3(n-\chi(M))$, it follows that 
\begin{equation}
\label{eq:delta}
\delta(G)\,\leq\, 6-\frac{6\chi(M)}{n}.
\end{equation}
If $\chi(M)>0$, then $\delta(G)\leq 5$. For $\chi(M)\leq 0$, we combine (\ref{eq:delta}) with $\delta(G)\leq n-1$
and obtain $\delta(G)\,\leq\, 6-\frac{6\chi(M)}{\delta(G)+1}$, which is equivalent to
\begin{equation}
\delta(G)\leq\Bigl\lfloor\tfrac{1}{2}(5+\sqrt{49-24\chi (M)})\Bigl\rfloor.
\end{equation}
In particular, we have $\delta(G)=q$ for any $q$-equivelar triangulation $K$ of a surface $M$,
and thus
\begin{equation}\label{eq:bound_for_q}
q\leq\Bigl\lfloor\tfrac{1}{2}(5+\sqrt{49-24\chi (M)})\Bigl\rfloor.
\end{equation}

Let $K$ be a $d$-covered triangulation of a surface $M$, $G=\mbox{\rm Sk}_1(K)$, and $H$ the subgraph of $G$ 
induced by the vertices of degree $d$. If $v$ is a vertex of $H$, then every edge of the link of $v$ in $G$ 
contains at least one vertex of $H$, since otherwise this edge would not be covered by a vertex of degree $d$. 
Hence, ${\rm deg}_{\,H}(v)\geq \frac{d}{2}$, from which it follows that 
\begin{equation}
\frac{d}{2}\leq \delta(H).
\end{equation}

If $K$ is a $3$-covered triangulation and $v$ a vertex of degree $3$, then the link of $v$ contains at least 
two vertices of degree $3$, which immediately implies that the third vertex of the link also
has degree~$3$, and thus $K$ is the boundary of the tetrahedron with $4$ vertices. 
If $K$ is a $4$-covered triangulation, it can be concluded that $K$ is the suspension
of a triangulated circle. For a characterization of $5$- and $6$-covered triangulations see \cite{NegamiNakamoto2001}.

Let $K$ be a $d$-covered triangulation of a surface $M$ and $w$ a vertex of degree ${\rm deg}_{\,G}(w)\neq d$.
Then all vertices of the link of $w$ are of degree $d$ and thus are vertices of $H$. Therefore, if we delete 
all vertices $w$ of degree ${\rm deg}_{\,G}(w)\neq d$ and the edges these vertices are contained in 
from the triangulation $K$, we obtain a regular cellular decomposition $C$ of $M$, with $H=\mbox{\rm Sk}_1(C)$. 
(A cell decomposition of a surface is \emph{regular} if there are no identifications on the boundary of any cell.)
Let $f^C=(f_0^C,f_1^C,f_2^C)$ be the face vector of $C$. For a $2$-face $P$ of $C$ let $p(P)$ be the number of vertices
of the polygon $P$, where $p(P)\geq 3$. Then, by the regularity of the cell complex~$C$, 
$2f_1^C=\sum_{P\in\{\mbox{\rm $2$-faces of $C$}\}}p(P)\geq 3f_2^C$.
Hence, $\chi(M)=f_0^C-f_1^C+f_2^C\leq f_0^C-\frac{1}{3}f_1^C$, or, equivalently, $f_1^C\leq 3(f_0^C-\chi(M))$.
It then follows that
\begin{equation}
\frac{d}{2}\leq\delta(H)=\frac{f_0^C\delta(H)}{f_0^C}\leq\frac{\sum_{\,v\in H\,}{\rm deg}_{\,H}(v)}{f_0^C}=\frac{2f_1^C}{f_0^C}\leq 6-\frac{6\chi(M)}{f_0^C},
\end{equation}
which, together with $\frac{d}{2}\leq \delta(H)\leq f_0^C-1$, gives 
\begin{equation}\label{eq:bound_for_d}
d\leq2\Bigl\lfloor\tfrac{1}{2}(5+\sqrt{49-24\chi (M)})\Bigl\rfloor
\end{equation}
if $\chi(M)\leq 0$, whereas $\frac{d}{2}\leq 5$ for $\chi(M)>0$.

\begin{thm}\label{thm:NegamiNakamotoI}
{\rm (Negami and Nakamoto \cite{NegamiNakamoto2001})}
Let $K$ be a $d$-covered triangulation of a surface $M$.
\begin{itemize}
\item[{\rm (a)}] If $\chi(M)>0$, then\, $d\leq 10$.
\item[{\rm (b)}] If $\chi(M)\leq 0$, then\, $d\leq 2\Bigl\lfloor\tfrac{1}{2}(5+\sqrt{49-24\chi (M)})\Bigl\rfloor$.
\end{itemize}
\end{thm}

If $d\geq 13$, then $7\leq\lceil \frac{d}{2}\rceil\leq\delta(H)\leq 6-\frac{6\chi(M)}{f_0^C}$,
or, equivalently $f_0^C\leq -6\chi(M)$. As there are only finitely many regular 
cell decompositions $C$ of a surface $M$ with $f_0^C\leq -6\chi(M)$,
we obtain:

\begin{thm}\label{thm:NegamiNakamotoII}
{\rm (Negami and Nakamoto \cite{NegamiNakamoto2001})}
Let $d\geq 13$ and $M$ be a closed surface. Then $M$ has only finitely many $d$-covered triangulations.
\end{thm}

\emph{Remark:} The proofs of the Theorems~\ref{thm:NegamiNakamotoI} and~\ref{thm:NegamiNakamotoII} given in \cite{NegamiNakamoto2001}
are somewhat vague and imprecise. Nevertheless, the statements of the two theorems are correct, as we verified above. 
In particular, we have for regular cell decompositions:

\begin{thm}
Let $C$ be a regular cell decomposition of a surface $M$ and let $H=\mbox{\rm Sk}_1(C)$.
\begin{itemize}
\item[{\rm (a)}] If $\chi(M)>0$, then\, $\delta(H)\leq 5$.
\item[{\rm (b)}] If $\chi(M)\leq 0$, then\, $\delta(H)\leq \Bigl\lfloor\tfrac{1}{2}(5+\sqrt{49-24\chi (M)})\Bigl\rfloor$.
\end{itemize}
\end{thm}

\begin{thm}
Let $M$ be a closed surface. Then $M$ has only finitely many regular cell decompositions~$C$ with $\delta(\mbox{\rm Sk}_1(C))\geq 7$.
\end{thm}

Since any $q$-equivelar triangulation yields, by subdivision, a $2q$-covered triangulation,
we see by the bounds (\ref{eq:bound_for_q}) and (\ref{eq:bound_for_d}),
that $q$-equivelar triangulations with
\begin{equation}
q=\Bigl\lfloor\tfrac{1}{2}(5+\sqrt{49-24\chi (M)})\Bigl\rfloor
\end{equation}
give $2q$-covered triangulation with 
\begin{equation}
d=2q=2\Bigl\lfloor\tfrac{1}{2}(5+\sqrt{49-24\chi (M)})\Bigl\rfloor,
\end{equation}
which means tightness in inequality (\ref{eq:bound_for_d}).

In the neighborly cases $n\equiv 0,1,3,4\,{\rm mod}\,6$, with $n\geq 4$,
there are, according to Theorem~\ref{thm:Jungerman_Ringel} and Equation~(\ref{eq:chi_equiv}), 
neighborly triangulations of orientable surfaces $M$ with
\begin{eqnarray}\label{eq:chi_or}
\chi(M) & = & \left\{ \begin{array}{l@{\hspace{2.5mm}}l@{\hspace{3mm}}l}
                          -24k^2+14k   & \mbox{\rm for} & n=12k,\,\,\, k\geq 1,\\
                          -24k^2+2k+2  & \mbox{\rm for} & n=12k+3,\,\,\, k\geq 1,\\ 
                          -24k^2-2k+2  & \mbox{\rm for} & n=12k+4,\,\,\, k\geq 0,\\
                          -24k^2-14k   & \mbox{\rm for} & n=12k+7,\,\,\, k\geq 0
                      \end{array}\right. 
\end{eqnarray}
and neighborly triangulations of non-orientable surfaces $M$ with
\begin{eqnarray}\label{eq:chi_non_or}
\chi(M) & = & \left\{ \begin{array}{l@{\hspace{2.5mm}}l@{\hspace{3mm}}l}
                          -\frac{3}{2}m^2+\frac{7}{2}m    & \mbox{\rm for} & n=3m,\,\, m\geq 2,\\[2mm]
                          -\frac{3}{2}m^2+\frac{5}{2}m+1  & \mbox{\rm for} & n=3m+1,\,\, m\geq 3.
                      \end{array}\right. 
\end{eqnarray}
The respective examples are equivelar with $q=n-1$ and yield $2q$-covered triangulations with tight bound~(\ref{eq:bound_for_d}).

\pagebreak

Apart from the series of neighborly triangulations of Ringel \cite{Ringel1955} and Jungerman and Ringel~\cite{JungermanRingel1980},
infinite classes and series of equivelar triangulation are known for all $q\geq 6$. In addition, there are some 
sporadic examples as well as enumerational results for small values of $n$ and $|\chi|$:

\begin{itemize}
\item The boundaries of the tetrahedron, the octahedron, and the icosahedron are the only equivelar triangulations 
      of the $2$-sphere $S^2$ with $\chi(S^2)=2$. In these cases, $(n,q)=(4,3)$, $(6,4)$, and $(12,5)$, respectively.
      The unique neighborly $6$-vertex triangulation of the real projective plane ${\mathbb R}{\bf P}^2$
      with $\chi({\mathbb R}{\bf P}^2)=1$ and $(n,q)=(6,5)$ is the only equivelar triangulation of ${\mathbb R}{\bf P}^2$.
      Together, these four examples are the only equivelar triangulations with $q\leq 5$.
\item There are infinitely many equivelar triangulations of the $2$-torus $T^2$ with $\chi(T^2)=0$ and $q=6$.
      A complete classification of these examples is due to Brehm and K\"uhnel \cite{BrehmKuehnel2008}.
      All of these examples have vertex-transitive cyclic symmetry \cite{DattaUpadhyay2005}. 

\begin{figure}
\begin{center}
\psfrag{0}{0}
\psfrag{1}{1}
\psfrag{2}{2}
\psfrag{3}{3}
\psfrag{4}{4}
\psfrag{5}{5}
\psfrag{n-2}{$n-2$}
\psfrag{n-1}{$n-1$}
\includegraphics[width=.675\linewidth]{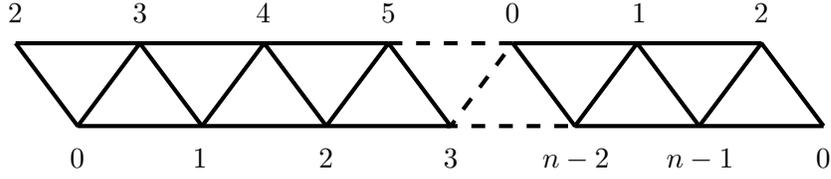}
\end{center}
\caption{The series of triangulated tori $T^2(n)$.}
\label{fig:T_2_n}
\end{figure}

      Let for $n\geq 7$ the cyclic group ${\mathbb Z}_n$ on the elements $0,1,\dots,n-1$ be generated
      by the permutation $(0,1,\dots,n-1)$ and let $T^2(n)$ be defined by the orbits of the triangles
      $[0,1,3]$ and $[0,2,3]$. Then $T^2(n)$, as depicted in Figure~\ref{fig:T_2_n}, is a cyclic series 
      of triangulated tori \cite{Altshuler1971},~\cite{KuehnelLassmann1996-bundle}. 
      For $n\geq 7$, the example $T^2(n)$ is realizable geometrically in the $2$-skeleton of the cyclic 
      $3$-polytope $C_3(n)$ \cite{Altshuler1971}.

\item Series of \emph{equivelar polyhedra} in $3$-space (i.e., geometric realizations in ${\mathbb R}^3$ of orientable equivelar surfaces) 
      were first described in \cite{McMullenSchulzWills1982}: There are equivelar polyhedra of the types~$\{3,q;g\}$ for 
      \begin{compactitem}
      \item[(a)] $\{3,7;g\}$ with\, $n=12(g-1)$\, for\, $g\geq 2$, 
      \item[(b)] $\{3,8;g\}$ with\, $n=6(g-1)$\, for\, $g\geq 4$,
      \item[(b)] $\{3,9;g\}$ with\, $n=4(g-1)$\, for\, $g=6,9,10$\, and\, $g\geq 12$, 
      \item[(d)] $\{3,12;g\}$ with\, $n=2(g-1)$\, for\, $g=73+66k$, $k\geq 0$,
      \end{compactitem} 
      with $g=1-\frac{1}{2}\chi$ the genus of a respective equivelar triangulation of
      vertex-degree $q$ and, as~a consequence of (\ref{eq:q}), with $n=\frac{6(2-2g)}{6-q}$ vertices.

\item McMullen, Schulz, and Wills \cite{McMullenSchulzWills1983} provided in 1983 a combinatorial construction method 
      for equivelar triangulations, which, in a second inductive step, allows a large family of equivelar polyhedra
      to be obtained.  
      In particular, for each $q\geq 7$ there are infinitely many combinatorially different equivelar polyhedra 
      of type $\{3,q\}$. Since for $g\geq 2$ the orientable surface $M_g$ of genus $g$ has (at most) 
      finitely many equivelar triangulations, it follows that for each $q\geq 7$ there are infinitely many equivelar polyhedra
      of distinct genus in the family.
 
      For $r\geq 3$, the vertex-minimal equivelar examples of type $\{3,2r\}$ in the constructed
      family have $n(r)=9\cdot2^{r-3}$ vertices and are of genus $g(r)=3(r-3)2^{r-4}+1$.
      By showing that for each $r\geq 3$ and all $0\leq g\leq g(r)$ there are equivelar polyhedra $M_{g}$ of genus $g$  
      with $f_0(M_{g})\leq n(r)$ vertices, McMullen, Schulz, and Wills established the remarkable result 
      that there are (equivelar) polyhedra of genus $g$ with $f_0=O(g/{\rm log}\,g)$ vertices.
      Asymptotically, this is the strongest known result for the realizability in ${\mathbb R}^3$ 
      of triangulated (orientable) surfaces of genus $g$, for which, by Heawood's bound (\ref{eq:Heawood}),
      $O(\sqrt{g})$ vertices are needed for a triangulation. A generalization of the McMullen-Schulz-Wills family 
      of equivelar surfaces is given in \cite{RoerigZiegler2009pre}.

\item Equivelar triangulations of the Klein bottle are of type $\{3,6\}$. The respective $1$-skeleta are $6$-regular graphs
      and were classified in \cite{Negami1984}; see also \cite{Fijavz2004} for additional comments and further references.
      Equivelar triangulations of the Klein bottle (and of the torus) with up to $15$ vertices were classified in 
      \cite{DattaUpadhyay2005}. An enumeration for up to $100$ vertices can be found in~\cite{SulankeLutz2009}.
      As observed in \cite{DattaUpadhyay2005},
      there is an $n$-vertex equivelar triangulation of the Klein bottle
      if and only if $n\geq 9$ is not prime.

\item Equivelar triangulations of the double torus have $12$ vertices and are of type $\{3,7\}$ by (\ref{eq:q}) and (\ref{eq:Heawood}):
      there are six such examples \cite{DattaUpadhyay2006}.

\item Equivelar triangulations with up to $11$ vertices were classified in \cite{DattaNilakantan2001}: there are $27$ such examples.
      By enumeration, there are $240\,914$ equivelar triangulations with $12$ vertices \cite{SulankeLutz2009}.

\item There are no equivelar triangulations of the non-orientable surface of genus~$3$ with $\chi=-1$.
      A complete enumeration of the equivelar triangulations 
      of the non-orientable surface of genus~$4$ with $\chi=-2$ ($28$ examples of type $(3,q;n)=(3,7;12)$), 
      of the non-orientable surface of genus $5$ with $\chi=-3$ (two examples of type $(3,8;9)$ and $1401$ examples of type $(3,7;18)$), 
      of the non-orientable surface of genus $6$ with $\chi=-4$ ($6500$ examples of type $(3,8;12)$ and $600\,946$ examples of type $(3,7;24)$),
      and of the orientable surface of genus $3$ with $\chi=-4$ ($24$ examples of type $(3,8;12)$ and $11\,277$ examples of type $(3,7;24)$)
      was given in \cite{SulankeLutz2009}.

\item The neighborly triangulations with up to $13$ vertices are: the boundary of the tetrahedron with $4$ vertices
      as the unique neighborly triangulation of $S^2$, the unique M\"obius torus \cite{Moebius1886} with $7$ vertices, 
      $59$ neighborly triangulations of the orientable surface of genus $6$ with $12$ vertices \cite{AltshulerBokowskiSchuchert1996}, 
      the unique vertex-minimal $6$-vertex triangulation of the real projective plane, two neighborly triangulations of the non-orientable surface
      of genus $5$ with $9$ vertices \cite{AltshulerBrehm1992}, $14$ neighborly triangulations of the non-orientable surface 
      of genus $7$ with $10$ vertices \cite{AltshulerBrehm1992}, $182\,200$ neighborly triangulations of the non-orientable 
      surface of genus $12$ with $12$ vertices \cite{EllinghamStephens2005}, \cite{SulankeLutz2009}, 
      and $243\,088\,286$ neighborly triangulations of the non-orientable surface of genus~$15$ with $13$ vertices \cite{EllinghamStephens2005}.
      
      The $4$-vertex triangulation of $S^2$ and M\"obius' torus are the only neighborly triangulations for which geometric realizations
      are known in ${\mathbb R}^3$ (for the Cs\'asz\'ar torus see \cite{Csaszar1949} and also \cite{BokowskiEggert1991}, \cite{Lutz2002b}). 
      None of the $59$ neighborly triangulations of the orientable surface of genus $6$ with $12$ vertices is realizable in ${\mathbb R}^3$ 
      \cite{BokowskiGuedes_de_Oliveira2000}, \cite{Schewe2008pre}, and it is assumed that neighborly triangulations
      of orientable surfaces of higher genus never are realizable. This leaves a wide gap to the realizability result
      of McMullen, Schulz, and Wills \cite{McMullenSchulzWills1983} that there are realizations of orientable surface
      of genus $g$ with $O(g/{\rm log}\,g)$ vertices.

\item A particular class of symmetric equivelar triangulations are vertex-transitive triangulations.
      An enumeration of vertex-transitive triangulations of surfaces with up to $15$ vertices
      was given in \cite{KoehlerLutz2005pre}, \cite{Lutz1999} and of vertex-transitive neighborly triangulations of surfaces
      with up to $22$ vertices in \cite{Lutz1999}.

\item Regular simplicial maps with a flag-transitive automorphism group are equivelar with highest possible symmetry.
      The book by Coxeter and Moser \cite{CoxeterMoser1957} gives a classical treatment. A~complete enumeration
      of all regular maps with Euler characteristic $\chi\geq -200$ was obtained by Conder \cite{Conder2009}.
      Dyck's regular map \cite{Dyck1880a}, \cite{Dyck1880b} and Klein's regular map \cite{Klein1879} are well-known examples 
      of regular simplicial maps for which geometric realization are known \cite{Bokowski1989}, \cite{BokowskiWills1988}, 
      \cite{Brehm1987a}, \cite{SchulteWills1985}, \cite{SchulteWills1986a}. 

\enlargethispage*{9.5mm}

\item There are two well-known infinite series of vertex-transitive triangulations of surfaces, the above series $T^2(n)$ of cyclic
      torus triangulations \cite{Altshuler1971},~\cite{KuehnelLassmann1996-bundle} 
      and Ringel's cyclic $7{\rm mod}12$ series of neighborly orientable triangulations \cite{Ringel1961}. 
      In fact, it was shown by Heffter \cite{Heffter1891} that cyclic neighborly triangulations of orientable surfaces 
      can only exist for $n\equiv 7{\rm mod}12$. The Altshuler series, the Ringel series, and the class of cyclic equivelar tori 
      \cite{BrehmKuehnel2008} are the only infinite series of vertex-transitive triangulations of surfaces that, so far, 
      have been described explicitely in the literature, while additional series of cyclic triangulations of surfaces
      of higher genus are implicitly contained in \cite{Brehm1990}, \cite{BrehmDattaNilakantan2002}, \cite{Datta2005}, and \cite{Jamet2001}.

\end{itemize}

Every equivelar triangulation of type $\{3,q\}$ gives rise to a $2q$-covered triangulation. 
Other sources for $d$-covered triangulations are:

\begin{itemize}

\item The boundary of the tetrahedron is the only $3$-covered triangulation. Any $4$-covered triangulation
      is the suspension of a triangulated circle. A characterization of $5$- and $6$-covered triangulations 
      is due to Negami and Nakamoto \cite{NegamiNakamoto2001}. For given $5\leq d\leq 12$, Katahira \cite{Katahira2004}
      proved that there are infinitely many distinct $d$-covered triangulations.

\item A \emph{map} on a surface $M$ is a decomposition of $M$ into a finite cell complex.
      A map is \emph{equivelar of type $\{p,q\}$} if $M$ is decomposed into $p$-gons such
      that every vertex has degree $q$; cf.~\cite{McMullenSchulzWills1982},~\cite{McMullenSchulzWills1983}. 
      A map is \emph{polyhedral} if the intersection of any two of its polygons is either empty, a common vertex, 
      or a common edge; see \cite{BrehmSchulte1997}, \cite{BrehmWills1993}.
      An \emph{equivelar polyhedral map} is a map, which is both equivelar and polyhedral. 

      If we place for each $p$-gon of a $\{p,q\}$-equivelar polyhedral map a new vertex inside the $p$-gon and connect 
      the new vertex with the $p$ vertices of the $p$-gon by an edge each, then the resulting triangulation is $2q$-covered. 
      
      Equivelar quadrangulations with few vertices are enumerated in \cite{LutzSulanke2009pre}, \cite{SulankeLutz2009}.
      McMullen, Schulz, and Wills~\cite{McMullenSchulzWills1983} gave infinite families of equivelar polyhedral maps 
      of the types $\{4,q\}$ and $\{p,4\}$ along with geometric realizations in ${\mathbb R}^3$
      of these examples. For geometric realizations of further equivelar polyhedral maps of these types 
      see \cite{BokowskiWills1988}, \cite{McMullenSchulteWills1988}, \cite{McMullenSchulzWills1982}, 
      \cite{SchulteWills1986b}, and \cite{RoerigZiegler2009pre}.

      It is not known whether there are geometric realizations of equivelar
      polyhedral maps of type $\{p,q\}$ for $p,q\geq 5$; cf.\ \cite{BrehmWills1993}. 
      Examples of equivelar polyhedral maps of type $\{5,5\}$ and of type $\{6,6\}$ were first
      given by Brehm~\cite{Brehm1990}. An infinite series of $\{k,k\}$-equivelar polyhedral maps 
      was constructed by Datta~\cite{Datta2005}.

\item A $\{p,q\}$-equivelar map is \emph{regular} if it has a flag-transitive automorphism group. 
      Regular polyhedral maps therefore provide highly symmetric examples of equivelar polyhedral maps;
      see \cite{Conder2009}, \cite{CoxeterMoser1957}.

\item If a $\{p,q\}$-equivelar map is not a polyhedral map, but at least gives a regular cell decomposition 
      of a surface, then this suffices to obtain a $2q$-covered triangulation of the surface by subdividing the $p$-gons.
      The Heffter surfaces \cite{Heffter1891} and the generalized Heffter surfaces 
      by Pfeifle and Ziegler \cite{PfeifleZiegler2002} provide such examples.

\end{itemize}

Infinite series of vertex-transitive triangulations of surfaces as well as of higher-dimensional manifolds 
are rare in the literature. The cyclic $7{\rm mod}12$ series of Ringel~\cite{Ringel1961} from 1961,
the cyclic series of triangulated tori of Altshuler~\cite{Altshuler1971} from 1971, and the class of all cyclic tori, 
as completely classified by Brehm and K\"uhnel \cite{BrehmKuehnel2008} in 2008,
are the only documented series of vertex-transitive triangulations of surfaces.
In higher dimensions, the boundaries of simplices, of cross-polytopes, and of cyclic polytopes
are classical examples of vertex-transitive triangulations of spheres.
A~cyclic series of centrally symmetric $3$-spheres is given in~\cite{Lutz2004apre},
series of wreath products are described in \cite{JoswigLutz2005}, vertex-transitive triangulations 
of higher-dimensional tori in \cite{KuehnelLassmann1988-dtori}, and series of cyclic triangulations 
of sphere bundles in \cite{Kuehnel1986a-series}, \cite{KuehnelLassmann1996-bundle}.

In the following section, we consider surfaces with Euler characteristic $0>\chi\geq -230$.
We show that in all non-neighborly cases with $q=\Bigl\lfloor\tfrac{1}{2}(5+\sqrt{49-24\chi (M)})\Bigl\rfloor$,
that is, with equality in~(\ref{eq:bound_for_q}), there indeed are equivelar triangulations
of the respective surfaces. By subdividing these examples, we obtain $2q$-covered triangulations 
with equality in~(\ref{eq:bound_for_d}).

In Section~\ref{sec:generalized_Ringel}, we generalize Ringel's cyclic $7{\rm mod}12$ series
to a two-parameter family $R_{k,n}$ of cyclic triangulations of the orientable surfaces 
of genus $g=kn+1$ with $n\geq 7+12k$ vertices and $q=6+12k$, $k\geq 0$. For $k=0$, $R_{0,n}$
coincides with Altshuler's cyclic series of tori, whereas $R_{k,7+12k}$ is Ringel's cyclic series
of neighborly triangulations with $n=7+12k$ vertices.
We show that, apart from the cyclic tori with $n\geq 8$ vertices, the two subseries $R_{k,7+12k+1}$ and $R_{k,7+12k+2}$, $k\geq 1$,
provide infinitely many examples of non-neighborly cyclic triangulations 
with $q=6+12k=\Bigl\lfloor\tfrac{1}{2}(5+\sqrt{49-24\chi})\Bigl\rfloor$,
i.e., with equality in~(\ref{eq:bound_for_q}). Again, by subdivision, these examples
yield $2q$-covered triangulations with equality in~(\ref{eq:bound_for_d}).

In addition to the generalized Ringel two-parameter family of cyclic triangulations, 
a generation scheme for a large class of series of cyclic triangulations 
of surfaces is presented in \cite{Lutz2010apre}. In particular,
infinite series of cyclic $q$-equivelar triangulations of orientable and 
non-orientable surfaces are given in \cite{Lutz2010apre} for each 
\begin{equation}
q = \left\{\begin{array}{lll}
                3k        & \mbox{\rm for} & k\geq 2,\\
                3k+1      & \mbox{\rm for} & k\geq 3.
               \end{array}\right.
\end{equation}

\section{\mathversion{bold}$d$-covered triangulations with $d=2\Bigl\lfloor\tfrac{1}{2}(5+\sqrt{49-24\chi (M)})\Bigl\rfloor$\mathversion{normal}}
\label{sec:small_examples}

In the orientable neighborly triangulated cases with $n\equiv 0,3,4,7\,{\rm mod}\,12$ vertices, $n\geq 4$,
and in the non-orientable neighborly triangulated cases with $n\equiv 0,1\,{\rm mod}\,3$ vertices, $n\geq 8$,
the Euler characteristic $\chi(M)$ of a respective orientable or non-orientable surface $M$ is determined,
as a consequence of (\ref{eq:chi_equiv}), by (\ref{eq:chi_or}) and (\ref{eq:chi_non_or}), respectively. 
In all these cases of neighborly triangulations, $n=\tfrac{1}{2}(7+\sqrt{49-24\chi (M)})$, $q=n-1$, 
and therefore $d=2q=2n-2=2\cdot \tfrac{1}{2}(5+\sqrt{49-24\chi (M)})$ for a derived $d$-covered triangulation.

In contrast, if $-\chi(M)\geq 5$ is prime, then Equation (\ref{eq:q}), $q=6-\frac{6\chi(M)}{n}$,
yields integer values for~$q$ if and only if $n=-\chi(M),-2\chi(M),-3\chi(M),-6\chi(M)$, with corresponding 
values $q=12,9,8,7$, respectively.

\begin{prop}
Let $M$ be a surface with $-\chi(M)\geq 5$ prime. If $M$ has an equivelar triangulation of type $\{3,q\}$,
then $q\in\{7,8,9,12\}$.
\end{prop}

\begin{table}
\centering
\defaultaddspace=0.15em
\caption{Existence of non-neighborly examples with equality in (\ref{eq:bound_for_q}) for $0>\chi\geq -127$.}\label{tbl:parameters}
\begin{tabular*}{\linewidth}{@{}r@{\extracolsep{28pt}}r@{\extracolsep{28pt}}r@{\extracolsep{\fill}}c@{\extracolsep{28pt}}c@{\extracolsep{\fill}}c@{\extracolsep{28pt}}c@{}}
\\\toprule
 \addlinespace
 \addlinespace
 \addlinespace
 \addlinespace
  $\chi$ &   $q$  &   $n$  &  Trans.\ or. &  Or.  &  Trans.\ non-or. &  Non-or. \\
\midrule
\\[-4mm]
 \addlinespace
 \addlinespace
 \addlinespace
 \addlinespace

   $-2$  &   $7$  &  $12$  &  yes \cite{KoehlerLutz2005pre} & yes \cite{SulankeLutz2009}    & yes \cite{KoehlerLutz2005pre} & yes \cite{SulankeLutz2009} \\
   $-4$  &   $8$  &  $12$  &  yes \cite{KoehlerLutz2005pre} & yes \cite{SulankeLutz2009}    & no  \cite{KoehlerLutz2005pre} & yes \cite{SulankeLutz2009} \\
   $-6$  &   $9$  &  $12$  &  yes \cite{KoehlerLutz2005pre} & yes \cite{SulankeLutz2009}    & yes \cite{KoehlerLutz2005pre} & yes \cite{SulankeLutz2009} \\
   $-7$  &   $9$  &  $14$  &  ---                           & ---                           & no  \cite{KoehlerLutz2005pre} & yes, E1                    \\
   $-8$  &  $10$  &  $12$  &  yes \cite{KoehlerLutz2005pre} & yes \cite{SulankeLutz2009}    & yes \cite{KoehlerLutz2005pre} & yes \cite{SulankeLutz2009} \\
  $-14$  &  $12$  &  $14$  &  yes \cite{KoehlerLutz2005pre} & yes \cite{KoehlerLutz2005pre} & yes \cite{KoehlerLutz2005pre} & yes \cite{KoehlerLutz2005pre}         \\
  $-15$  &  $12$  &  $15$  &  ---                           & ---                           & yes \cite{KoehlerLutz2005pre} & yes \cite{KoehlerLutz2005pre}         \\
  $-16$  &  $12$  &  $16$  &  yes, E2                       & yes, E2                       & yes, E3                       & yes, E3  \\
  $-27$  &  $15$  &  $18$  &  ---                           & ---                           & no                            & yes, E4  \\
  $-30$  &  $16$  &  $18$  &  yes, E5                       & yes, E5                       & yes, E6                       & yes, E6  \\
  $-40$  &  $18$  &  $20$  &  yes, E7                       & yes, E7                       & yes, E8                       & yes, E8  \\
  $-42$  &  $18$  &  $21$  &  yes, E9                       & yes, E9                       & yes, E10                      & yes, E10 \\
  $-60$  &  $21$  &  $24$  &  yes, E11                      & yes, E11                      & ?                             & yes, E12 \\
  $-64$  &  $22$  &  $24$  &  ?                             & yes, E13                      & yes, E14                      & yes, E14 \\
  $-78$  &  $24$  &  $26$  &  yes, E15                      & yes, E15                      & ?                             & yes, E16 \\
  $-81$  &  $24$  &  $27$  &  ---                           & ---                           & yes, E17                      & yes, E17 \\
 $-105$  &  $27$  &  $30$  &  ---                           & ---                           & yes, E18                      & yes, E18 \\
 $-110$  &  $28$  &  $30$  &  ?                             & yes, E19                      & yes, E20                      & yes, E20 \\

 \addlinespace

 \addlinespace
 \addlinespace
 \addlinespace
 \addlinespace
 \bottomrule
\end{tabular*}
\end{table}

In \cite{Upadhyay2008pre}, the fourth author raised the question, whether there are cases, other than the neighborly ones,
where $d=2\lfloor\tfrac{1}{2}(5+\sqrt{49-24\chi (M)})\rfloor$ can be achieved by subdivisions
of $q$-equivelar triangulations with equality in (\ref{eq:bound_for_q}), that is, 
with $q=\lfloor\tfrac{1}{2}(5+\sqrt{49-24\chi (M)})\rfloor$.

For $\chi=0$, all equivelar triangulations of the torus and the Klein bottle are of type $\{3,6\}$,
and thus yield $12$-covered triangulations. So let $\chi<0$. Then for a $q$-equivelar triangulation
with $q=\lfloor\tfrac{1}{2}(5+\sqrt{49-24\chi (M)})\rfloor$ we have by (\ref{eq:q}) that
\begin{equation}
n\,\,=\,\,-\frac{6\chi(M)}{q-6}\,\,=\,\,\frac{-6\chi(M)}{\lfloor\tfrac{1}{2}(5+\sqrt{49-24\chi (M)})\rfloor-6},
\end{equation}
which needs to be an integer.

\enlargethispage*{11mm}

Table~\ref{tbl:parameters} lists for $0>\chi\geq -127$ all those triples $(\chi,q,n)$ with 
values $n\geq \Bigl\lceil\tfrac{1}{2}(7+\sqrt{49-24\chi})\Bigl\rceil$,
for which $q=\lfloor\tfrac{1}{2}(5+\sqrt{49-24\chi})\rfloor<n-1$
and $n=\frac{-6\chi}{\lfloor\tfrac{1}{2}(5+\sqrt{49-24\chi})\rfloor-6}$ is an integer. 
These are, apart from the neighborly cases, the potential cases where 
$q$-equivelar triangulations with equality $q=\lfloor\tfrac{1}{2}(5+\sqrt{49-24\chi})\rfloor$ can occur.
Columns four and six of Table~\ref{tbl:parameters} indicate whether vertex-transitive orientable and non-orientable
examples are known, whereas columns five and seven provide information on orientable and non-orientable equivelar
triangulations, respectively. The vertex-transitive examples E2--E3, E5--E11, E14--E15, E17--E18, and E20 are listed in Table~\ref{tbl:examples},
the non-transitive examples E1, E4, E12, E13, E16, and E19 are listed at the end of this section.

\begin{table}
\centering
\defaultaddspace=0.15em
\caption{Examples with vertex-transitive symmetry.}\label{tbl:examples}
\begin{tabular*}{\linewidth}{@{}r@{\extracolsep{\fill}}r@{\extracolsep{\fill}}l@{}}
\\\toprule
 \addlinespace
 \addlinespace
 \addlinespace
 \addlinespace
  Example  & Group: $(n,i)$  &  Orbit generating triangles \\
\midrule
\\[-4mm]
 \addlinespace
 \addlinespace
 \addlinespace
 \addlinespace

   E2  &  (16,189)  &  $[1,3,5]_{16}$, $[1,3,7]_{48}$  \\ 
   E3  &  (16,28)   &  $[1,3,6]_{32}$, $[1,3,9]_{32}$  \\ 
   E5  &  (18,3)    &  $[1,2,5]_{18}$, $[1,3,8]_{18}$, $[1,4,9]_{18}$, $[1,4,17]_{6}$, $[1,5,12]_{18}$, $[1,6,11]_{18}$  \\ 
   E6  &  (18,3)    &  $[1,2,3]_{18}$, $[1,4,6]_{18}$, $[1,5,11]_{18}$, $[1,5,15]_{6}$, $[1,6,16]_{18}$, $[1,7,13]_{18}$  \\ 
   E7  &  (20,35)   &  $[1,3,5]_{120}$ \\ 
   E8  &  (20,2)    &  $[1,3,5]_{20}$, $[1,3,7]_{20}$, $[1,6,13]_{20}$, $[1,6,15]_{20}$, $[1,7,13]_{20}$, $[1,10,19]_{20}$  \\ 
   E9  &  (21,7)    &  $[1,4,10]_{21}$, $[1,4,11]_{63}$, $[1,5,12]_{21}$, $[1,6,16]_{21}$  \\ 
   E10 &  (21,7)    &  $[1,4,11]_{63}$, $[1,4,13]_{63}$  \\ 
   E11 &  (24,84)   &  $[1,4,7]_{24}$, $[1,4,14]_{72}$, $[1,5,10]_{72}$  \\ 
   E14 &  (24,21)   &  $[1,3,8]_{48}$, $[1,3,23]_{16}$, $[1,4,13]_{48}$, $[1,4,15]_{48}$, $[1,5,16]_{16}$  \\ 
   E15 &  (26,5)    &  $[1,3,5]_{78}$, $[1,3,13]_{78}$, $[1,4,10]_{26}$, $[1,6,18]_{26}$  \\ 
   E17 &  (27,23)   &  $[1,4,10]_{81}$, $[1,4,11]_{81}$, $[1,8,18]_{27}$, $[1,9,17]_{27}$   \\ 
   E18 &  (30,9)    &  $[1,2,21]_{30}$, $[1,3,8]_{60}$, $[1,4,9]_{60}$, $[1,5,14]_{60}$, $[1,9,28]_{60}$  \\ 
   E20 &  (30,9)    &  $[1,3,8]_{60}$, $[1,4,9]_{60}$, $[1,5,14]_{60}$, $[1,9,30]_{30}$, $[1,10,30]_{10}$, $[1,15,17]_{60}$  \\ 
   
 \addlinespace

 \addlinespace
 \addlinespace
 \addlinespace
 \addlinespace
 \bottomrule
\end{tabular*}
\end{table}

\begin{thm}
There are exactly $42$ (orientable and non-orientable) surfaces with $0>\chi\geq -230$ 
for which there are $q$-equivelar non-neighborly triangulations with equality in~(\ref{eq:bound_for_q}) 
and derived $2q$-covered triangulations with equality in~(\ref{eq:bound_for_d}).
\end{thm}

\textbf{Proof:} A complete list of equivelar triangulations with up to $12$ vertices is given in \cite{SulankeLutz2009},
whereas vertex-transitive triangulations with up to $15$ vertices can be found in \cite{KoehlerLutz2005pre} and \cite{Lutz1999}.
For~\mbox{$16\leq n\leq 30$}, we used the GAP program MANIFOLD\_VT \cite{Lutz_MANIFOLD_VT} to search for vertex-transitive triangulations.

In most of the cases of Table~\ref{tbl:parameters}, corresponding vertex-transitive 
examples were found. However, we know from \cite{KoehlerLutz2005pre}, \cite{Lutz1999} that there are no vertex-transitive 
triangulations  of the non-orientable surface of Euler characteristic $\chi=-4$ with $n=12$ vertices 
and no vertex-transitive triangulations of the surface of Euler characteristic $\chi=-7$ with \mbox{$n=14$}. 
It also turned out that there are no vertex-transitive triangulations of the surface of Euler characteristic 
$\chi=-27$ with $n=18$. Nevertheless, there are equivelar triangulations in these three cases, 
$500$ non-orientable equivelar examples in the case $(\chi,q,n)=(-4,8,12)$, $11\,300\,559$~equivelar examples 
in the case $(-7,9,14)$ \cite{Lutz_PAGE}, and at least $412$ equivelar examples in the case $(-27,15,18)$, 
as were found with the program \texttt{lextri} \cite{Sulanke2009}. 
In the range $0>\chi\geq -127$ there remain four undecided cases where we do not know whether 
vertex-transitive triangulations exist. For the search for vertex-transitive triangulations 
we made use of the classification of transitive permutations groups of small degree~$n$. 
A~complete list of these groups is available for $n\leq 30$ via the package GAP~\cite{GAP4}, 
imposing a restriction $\chi\geq -127$ for the search (since the next non-neighborly candidate case with $\chi=-128$ 
requires $n=32$ vertices). For $n=24$, $26$, and~$30$, we were only able to process all transitive 
permutations groups of order larger than $n$ (as well as some of order $n$) with MANIFOLD\_VT . Thus, in the open cases, 
there might be vertex-transitive triangulations corresponding to transitive permutations groups 
of degree $n$ and order $n$. In these cases also \texttt{lextri} was too slow to produce examples,
but we were able to find respective examples by using bistellar flips (with pruning the search
according to the fact that before the final flip of the search exactly four vertices have the wrong degree).
The vertex-transitive examples we found are listed in Table~\ref{tbl:examples},
with $(n,i)$ the $i$th transitive permutation group of degree $n$ in the GAP catalog,
and the size of an orbit indicated in subscript.
The non-transitive examples E1, E4, E12, E13, E16, and E19 are listed at the end of this section.

In the range $31\leq n\leq 40$, that is, $-128\geq\chi\geq -230$, there are six further triples
$(-128,30,32)$, $(-132,30,33)$, $(-162,33,36)$, $(-168,34,36)$, $(-190,36,38)$, and $(-195,36,39)$.
We found orientable examples in the first five cases and non-orientable examples in all six cases.~\hfill$\Box$

\medskip

The non-orientable surface of Euler characteristic $\chi=-3$ is the smallest case, where for an admissible triple, 
$(\chi,q,n)=(-3,8,9)$, there are no vertex-transitive triangulations, whereas two equivelar triangulations 
exist; see \cite{KoehlerLutz2005pre}, \cite{Lutz1999}, and \cite{SulankeLutz2009}. 
Further cases without vertex-transitive, but with equivelar triangulations
are the non-orientable cases $(-4,8,12)$, $(-7,9,14)$, and $(-27,15,18)$.

\begin{conj}\label{conj:triples}
Let $\chi<0$, $n\geq\Bigl\lceil\tfrac{1}{2}(7+\sqrt{49-24\chi})\Bigl\rceil$, and $q=6-\frac{6\chi(M)}{n}$ an integer. 
For every admissible triple $(\chi,q,n)$, the non-orientable surface of Euler characteristic $\chi$ 
and, if $\chi$ is even, the orientable surface of Euler characteristic $\chi$ have 
$q$-equivelar triangulations with $n$ vertices.
\end{conj}

\enlargethispage*{6mm}

Conjecture~\ref{conj:triples} holds for orientable surfaces in the cases $q=6+12k$, $k\geq 1$; see Section~\ref{sec:generalized_Ringel}.
Moreover, we found respective examples for all admissible triples $(\chi,q,n)$
with $n\leq 40$ with the help of bistellar flips.

\begin{thm}
For each admissible triple $(\chi,q,n)$ with $n\leq 40$ there is a $q$-equivelar triangulation
with $n$ vertices of the non-orientable surface of Euler characteristic $\chi$ and, if $\chi$ is even, 
then also of the orientable surface of Euler characteristic $\chi$.
\end{thm}

For $\chi>0$, there is exactly one equivelar, vertex-transitive triangulation
of $S^2$ for $(n,q)=(4,3)$, $(6,4)$, and $(12,5)$ each, as well of ${\mathbb R}{\bf P}^2$ for $(n,q)=(6,5)$.
The $2$-torus has equivelar, vertex-transitive triangulations (with $q=6$) for all $n\geq 7$\, \cite{Altshuler1971}, \cite{BrehmKuehnel2008},
whereas the Klein bottle has equivelar triangulations (with $q=6$) for all non-prime $n\geq 9$\, \cite{DattaUpadhyay2005}.

\noindent
Example E1:\\[2mm]
{\scriptsize
\begin{tabular*}{\linewidth}{@{}l@{\extracolsep{\fill}}l@{\extracolsep{\fill}}l@{\extracolsep{\fill}}l@{\extracolsep{\fill}}l@{\extracolsep{\fill}}l@{\extracolsep{\fill}}l@{\extracolsep{\fill}}l@{\extracolsep{\fill}}l@{\extracolsep{\fill}}l@{}}
     123  &  124  &  135  &  146  &  157  &  168  &  179  &  18\,10  &  19\,10  &  236  \\
     245  &  257  &  268  &  27\,10  &  289  &  29\,10  &  347  &  348  &  35\,11  &  36\,12  \\
     37\,11  &  38\,12  &  45\,13  &  46\,14  &  47\,14  &  48\,13  &  5\,11\,12  &  5\,12\,14  &  5\,13\,14  &  6\,11\,12  \\
     6\,11\,13  &  6\,13\,14  &  79\,11  &  7\,10\,14  &  89\,13  &  8\,10\,12  &  9\,11\,14  &  9\,12\,13  &  9\,12\,14  &  10\,11\,13  \\
     10\,11\,14  &  10\,12\,13.
\end{tabular*}
}

\vspace{2mm}

\noindent
Example E4:\\[2mm]
{\scriptsize
\begin{tabular*}{\linewidth}{@{}l@{\extracolsep{\fill}}l@{\extracolsep{\fill}}l@{\extracolsep{\fill}}l@{\extracolsep{\fill}}l@{\extracolsep{\fill}}l@{\extracolsep{\fill}}l@{\extracolsep{\fill}}l@{\extracolsep{\fill}}l@{\extracolsep{\fill}}l@{}}
     12\,12  &  12\,17  &  134  &  13\,18  &  145  &  159  &  167  &  16\,14  &  17\,12  &  18\,11  \\
     18\,17  &  19\,15  &  1\,11\,15  &  1\,14\,16  &  1\,16\,18  &  236  &  23\,10  &  257  &  25\,14  &  26\,13  \\
     27\,16  &  289  &  28\,14  &  29\,17  &  2\,10\,11  &  2\,11\,18  &  2\,12\,13  &  2\,16\,18  &  34\,15  &  356  \\
     35\,17  &  38\,12  &  38\,14  &  39\,13  &  39\,18  &  3\,10\,13  &  3\,11\,12  &  3\,11\,15  &  3\,14\,17  &  45\,12  \\
     47\,14  &  47\,17  &  48\,15  &  48\,17  &  49\,10  &  49\,11  &  4\,10\,18  &  4\,11\,13  &  4\,12\,16  &  4\,13\,16  \\
     4\,14\,18  &  56\,10  &  57\,15  &  58\,12  &  58\,16  &  59\,10  &  5\,14\,18  &  5\,15\,16  &  5\,17\,18  &  67\,11  \\
     69\,12  &  69\,15  &  6\,10\,11  &  6\,12\,14  &  6\,13\,18  &  6\,15\,16  &  6\,16\,17  &  6\,17\,18  &  79\,12  &  79\,18  \\
     7\,10\,14  &  7\,10\,15  &  7\,11\,16  &  7\,13\,17  &  7\,13\,18  &  89\,13  &  8\,10\,16  &  8\,10\,18  &  8\,11\,18  &  8\,13\,15  \\
     9\,11\,17  &  10\,12\,15  &  10\,12\,16  &  10\,13\,14  &  11\,12\,13  &  11\,16\,17  &  12\,14\,15  &  13\,14\,16  &  13\,15\,17  &  14\,15\,17.
\end{tabular*}
}

\vspace{2mm}

\noindent
Example E12:\\[2mm]
{\scriptsize
\begin{tabular*}{\linewidth}{@{}l@{\extracolsep{\fill}}l@{\extracolsep{\fill}}l@{\extracolsep{\fill}}l@{\extracolsep{\fill}}l@{\extracolsep{\fill}}l@{\extracolsep{\fill}}l@{\extracolsep{\fill}}l@{\extracolsep{\fill}}l@{\extracolsep{\fill}}l@{}}
     12\,11  &  12\,20  &  134  &  13\,10  &  14\,23  &  15\,22  &  15\,24  &  16\,17  &  16\,21  &  17\,11  \\
     17\,12  &  18\,16  &  18\,17  &  19\,22  &  19\,23  &  1\,10\,16  &  1\,12\,18  &  1\,15\,19  &  1\,15\,21  &  1\,18\,24  \\
     1\,19\,20  &  23\,14  &  23\,15  &  247  &  24\,22  &  25\,16  &  25\,18  &  27\,23  &  28\,15  &  28\,17  \\
     29\,16  &  29\,23  &  2\,10\,17  &  2\,10\,21  &  2\,11\,12  &  2\,12\,21  &  2\,13\,14  &  2\,13\,19  &  2\,18\,19  &  2\,20\,22  \\
     34\,14  &  35\,20  &  35\,24  &  36\,16  &  36\,23  &  379  &  37\,24  &  38\,12  &  38\,22  &  39\,12  \\
     3\,10\,18  &  3\,11\,21  &  3\,11\,23  &  3\,15\,16  &  3\,18\,20  &  3\,19\,21  &  3\,19\,22  &  457  &  45\,11  &  46\,13  \\
     46\,18  &  48\,12  &  48\,21  &  49\,11  &  49\,20  &  4\,12\,23  &  4\,13\,22  &  4\,14\,17  &  4\,15\,18  &  4\,15\,24  \\
     4\,16\,17  &  4\,16\,20  &  4\,21\,24  &  56\,11  &  56\,15  &  57\,22  &  589  &  58\,21  &  59\,13  &  5\,10\,18  \\
     5\,10\,20  &  5\,12\,13  &  5\,12\,15  &  5\,14\,16  &  5\,14\,17  &  5\,17\,21  &  67\,21  &  67\,22  &  69\,18  &  69\,19  \\
     6\,10\,12  &  6\,10\,23  &  6\,11\,20  &  6\,12\,14  &  6\,13\,17  &  6\,14\,22  &  6\,15\,19  &  6\,16\,24  &  6\,20\,24  &  78\,18  \\
     78\,24  &  79\,14  &  7\,10\,11  &  7\,10\,13  &  7\,12\,17  &  7\,13\,19  &  7\,14\,21  &  7\,17\,18  &  7\,19\,20  &  7\,20\,23  \\
     89\,19  &  8\,11\,13  &  8\,11\,16  &  8\,13\,18  &  8\,14\,19  &  8\,14\,24  &  8\,15\,20  &  8\,20\,23  &  8\,22\,23  &  9\,10\,12  \\
     9\,10\,16  &  9\,11\,17  &  9\,13\,24  &  9\,14\,24  &  9\,17\,22  &  9\,18\,20  &  10\,11\,19  &  10\,13\,20  &  10\,14\,15  &  10\,14\,23  \\
     10\,15\,17  &  10\,19\,24  &  10\,21\,22  &  10\,22\,24  &  11\,12\,16  &  11\,13\,15  &  11\,15\,23  &  11\,17\,24  &  11\,18\,19  &  11\,18\,24  \\
     11\,20\,21  &  12\,13\,21  &  12\,14\,18  &  12\,15\,22  &  12\,16\,17  &  12\,19\,22  &  12\,19\,23  &  13\,14\,22  &  13\,15\,16  &  13\,16\,23  \\
     13\,17\,23  &  13\,18\,21  &  13\,20\,24  &  14\,15\,20  &  14\,16\,19  &  14\,18\,23  &  14\,20\,21  &  15\,17\,22  &  15\,18\,21  &  15\,23\,24  \\
     16\,19\,23  &  16\,20\,22  &  16\,21\,22  &  16\,21\,24  &  17\,18\,23  &  17\,19\,21  &  17\,19\,24  &  22\,23\,24.
\end{tabular*}
}

\vspace{2mm}

\noindent
Example E13:\\[2mm]
{\scriptsize
\begin{tabular*}{\linewidth}{@{}l@{\extracolsep{\fill}}l@{\extracolsep{\fill}}l@{\extracolsep{\fill}}l@{\extracolsep{\fill}}l@{\extracolsep{\fill}}l@{\extracolsep{\fill}}l@{\extracolsep{\fill}}l@{\extracolsep{\fill}}l@{\extracolsep{\fill}}l@{}}
     124  &  12\,20  &  13\,13  &  13\,17  &  14\,14  &  16\,11  &  16\,22  &  178  &  17\,19  &  18\,10  \\
     19\,15  &  19\,24  &  1\,10\,12  &  1\,11\,21  &  1\,12\,21  &  1\,13\,16  &  1\,14\,16  &  1\,15\,22  &  1\,17\,23  &  1\,18\,23  \\
     1\,18\,24  &  1\,19\,20  &  235  &  23\,20  &  24\,12  &  25\,18  &  269  &  26\,12  &  278  &  27\,16  \\
     289  &  2\,10\,13  &  2\,10\,15  &  2\,11\,14  &  2\,11\,23  &  2\,13\,24  &  2\,14\,24  &  2\,15\,19  &  2\,16\,22  &  2\,17\,22  \\
     2\,17\,23  &  2\,18\,19  &  34\,16  &  34\,21  &  35\,24  &  36\,19  &  36\,21  &  37\,13  &  37\,15  &  389  \\
     38\,24  &  39\,11  &  3\,10\,11  &  3\,10\,23  &  3\,12\,18  &  3\,12\,23  &  3\,15\,19  &  3\,16\,22  &  3\,17\,20  &  3\,18\,22  \\
     45\,13  &  45\,17  &  46\,15  &  46\,22  &  479  &  47\,24  &  48\,10  &  48\,12  &  49\,16  &  4\,10\,20  \\
     4\,13\,23  &  4\,14\,18  &  4\,15\,17  &  4\,18\,19  &  4\,19\,22  &  4\,20\,24  &  4\,21\,23  &  56\,17  &  56\,20  &  579  \\
     57\,14  &  58\,20  &  58\,21  &  59\,22  &  5\,10\,15  &  5\,10\,16  &  5\,11\,12  &  5\,11\,18  &  5\,12\,19  &  5\,13\,19  \\
     5\,14\,21  &  5\,15\,23  &  5\,16\,24  &  5\,22\,23  &  68\,17  &  68\,24  &  69\,21  &  6\,10\,12  &  6\,10\,18  &  6\,11\,13  \\
     6\,13\,14  &  6\,14\,23  &  6\,15\,18  &  6\,16\,20  &  6\,16\,23  &  6\,19\,24  &  7\,10\,17  &  7\,10\,18  &  7\,11\,12  &  7\,11\,21  \\
     7\,12\,20  &  7\,13\,23  &  7\,14\,20  &  7\,15\,22  &  7\,16\,17  &  7\,18\,21  &  7\,19\,23  &  7\,22\,24  &  8\,11\,16  &  8\,11\,20  \\
     8\,12\,19  &  8\,13\,14  &  8\,13\,22  &  8\,14\,17  &  8\,16\,23  &  8\,18\,21  &  8\,18\,22  &  8\,19\,23  &  9\,10\,14  &  9\,10\,16  \\
     9\,11\,14  &  9\,12\,17  &  9\,12\,18  &  9\,13\,19  &  9\,13\,21  &  9\,15\,20  &  9\,17\,22  &  9\,18\,20  &  9\,19\,24  &  10\,11\,24  \\
     10\,13\,17  &  10\,14\,22  &  10\,20\,21  &  10\,21\,23  &  10\,22\,24  &  11\,13\,20  &  11\,15\,16  &  11\,15\,17  &  11\,17\,19  &  11\,18\,24  \\
     11\,19\,22  &  11\,22\,23  &  12\,14\,15  &  12\,14\,22  &  12\,15\,24  &  12\,16\,17  &  12\,16\,20  &  12\,21\,22  &  12\,23\,24  &  13\,15\,18  \\
     13\,15\,20  &  13\,16\,24  &  13\,17\,18  &  13\,21\,22  &  14\,15\,23  &  14\,16\,19  &  14\,17\,18  &  14\,19\,20  &  14\,21\,24  &  15\,16\,21  \\
     15\,21\,24  &  16\,19\,21  &  17\,19\,21  &  17\,20\,21  &  18\,20\,23  &  20\,23\,24.
\end{tabular*}
}

\vspace{3mm}

\noindent
Example E16:\\[3mm]
{\scriptsize
\begin{tabular*}{\linewidth}{@{}l@{\extracolsep{\fill}}l@{\extracolsep{\fill}}l@{\extracolsep{\fill}}l@{\extracolsep{\fill}}l@{\extracolsep{\fill}}l@{\extracolsep{\fill}}l@{\extracolsep{\fill}}l@{\extracolsep{\fill}}l@{\extracolsep{\fill}}l@{}}
     12\,15  &  12\,23  &  13\,15  &  13\,17  &  14\,10  &  14\,18  &  15\,17  &  15\,21  &  167  &  16\,22  \\
     17\,12  &  18\,13  &  18\,16  &  19\,13  &  19\,20  &  1\,10\,26  &  1\,11\,18  &  1\,11\,24  &  1\,12\,14  &  1\,14\,24  \\
     1\,16\,19  &  1\,19\,21  &  1\,20\,26  &  1\,22\,23  &  23\,19  &  23\,23  &  247  &  249  &  257  &  25\,11  \\
     26\,12  &  26\,24  &  28\,14  &  28\,26  &  29\,20  &  2\,10\,14  &  2\,10\,17  &  2\,11\,13  &  2\,12\,19  &  2\,13\,17  \\
     2\,15\,20  &  2\,16\,24  &  2\,16\,25  &  2\,21\,22  &  2\,21\,25  &  2\,22\,26  &  348  &  34\,18  &  358  &  35\,21  \\
     36\,16  &  36\,19  &  37\,13  &  37\,18  &  39\,14  &  39\,24  &  3\,10\,16  &  3\,10\,26  &  3\,11\,20  &  3\,11\,26  \\
     3\,12\,15  &  3\,12\,25  &  3\,13\,20  &  3\,14\,25  &  3\,17\,23  &  3\,21\,24  &  456  &  45\,13  &  469  &  47\,26  \\
     48\,20  &  4\,10\,23  &  4\,11\,15  &  4\,11\,24  &  4\,13\,16  &  4\,14\,19  &  4\,14\,26  &  4\,15\,22  &  4\,16\,20  &  4\,17\,24  \\
     4\,17\,25  &  4\,19\,21  &  4\,21\,22  &  4\,23\,25  &  56\,24  &  57\,20  &  58\,12  &  59\,19  &  59\,23  &  5\,10\,22  \\
     5\,10\,25  &  5\,11\,19  &  5\,12\,14  &  5\,13\,14  &  5\,16\,17  &  5\,16\,26  &  5\,18\,20  &  5\,18\,25  &  5\,22\,24  &  5\,23\,26  \\
     67\,26  &  68\,11  &  68\,15  &  69\,21  &  6\,11\,19  &  6\,12\,20  &  6\,13\,14  &  6\,13\,18  &  6\,14\,25  &  6\,15\,23  \\
     6\,16\,26  &  6\,17\,18  &  6\,17\,21  &  6\,20\,22  &  6\,23\,25  &  78\,10  &  78\,18  &  79\,12  &  79\,19  &  7\,10\,22  \\
     7\,13\,15  &  7\,14\,15  &  7\,14\,17  &  7\,16\,20  &  7\,16\,21  &  7\,17\,21  &  7\,19\,23  &  7\,22\,25  &  7\,23\,24  &  7\,24\,25  \\
     8\,10\,21  &  8\,11\,21  &  8\,12\,17  &  8\,13\,22  &  8\,14\,19  &  8\,15\,25  &  8\,16\,17  &  8\,18\,22  &  8\,19\,25  &  8\,20\,23  \\
     8\,23\,24  &  8\,24\,26  &  9\,10\,13  &  9\,10\,17  &  9\,11\,16  &  9\,11\,23  &  9\,12\,22  &  9\,14\,18  &  9\,15\,16  &  9\,15\,25  \\
     9\,17\,25  &  9\,18\,26  &  9\,21\,26  &  9\,22\,24  &  10\,11\,12  &  10\,11\,14  &  10\,12\,25  &  10\,13\,23  &  10\,15\,16  &  10\,15\,19  \\
     10\,18\,20  &  10\,18\,24  &  10\,19\,24  &  10\,20\,21  &  11\,12\,16  &  11\,13\,18  &  11\,14\,15  &  11\,17\,22  &  11\,17\,23  &  11\,20\,25  \\
     11\,21\,26  &  11\,22\,25  &  12\,13\,20  &  12\,13\,24  &  12\,15\,26  &  12\,16\,22  &  12\,17\,18  &  12\,18\,23  &  12\,19\,26  &  12\,21\,23  \\
     12\,21\,24  &  13\,15\,21  &  13\,16\,25  &  13\,17\,22  &  13\,21\,23  &  13\,24\,26  &  13\,25\,26  &  14\,16\,22  &  14\,16\,24  &  14\,17\,20  \\
     14\,18\,23  &  14\,20\,23  &  14\,22\,26  &  15\,17\,19  &  15\,17\,24  &  15\,18\,21  &  15\,18\,24  &  15\,20\,22  &  15\,23\,26  &  16\,18\,19  \\
     16\,18\,21  &  17\,19\,20  &  18\,19\,22  &  18\,25\,26  &  19\,20\,26  &  19\,22\,23  &  19\,24\,25  &  20\,21\,25.
\end{tabular*}
}

\vspace{3mm}

\noindent
Example E19:\\[3mm]
{\scriptsize
\begin{tabular*}{\linewidth}{@{}l@{\extracolsep{\fill}}l@{\extracolsep{\fill}}l@{\extracolsep{\fill}}l@{\extracolsep{\fill}}l@{\extracolsep{\fill}}l@{\extracolsep{\fill}}l@{\extracolsep{\fill}}l@{\extracolsep{\fill}}l@{\extracolsep{\fill}}l@{}}
     125  &  12\,15  &  13\,16  &  13\,25  &  14\,10  &  14\,11  &  15\,13  &  16\,11  &  16\,20  &  17\,15  \\
     17\,17  &  18\,20  &  18\,29  &  19\,17  &  19\,23  &  1\,10\,12  &  1\,12\,22  &  1\,13\,30  &  1\,14\,23  &  1\,14\,27  \\
     1\,16\,27  &  1\,19\,28  &  1\,19\,29  &  1\,21\,22  &  1\,21\,26  &  1\,24\,25  &  1\,24\,28  &  1\,26\,30  &  23\,20  &  23\,30  \\
     249  &  24\,21  &  258  &  26\,22  &  26\,25  &  27\,19  &  27\,29  &  28\,14  &  29\,16  &  2\,10\,12  \\
     2\,10\,13  &  2\,11\,23  &  2\,11\,30  &  2\,12\,28  &  2\,13\,15  &  2\,14\,17  &  2\,16\,22  &  2\,17\,25  &  2\,18\,24  &  2\,18\,29  \\
     2\,19\,28  &  2\,20\,24  &  2\,21\,27  &  2\,23\,27  &  346  &  349  &  35\,16  &  35\,24  &  36\,29  &  37\,14  \\
     37\,24  &  39\,21  &  3\,10\,19  &  3\,10\,28  &  3\,11\,14  &  3\,11\,27  &  3\,12\,17  &  3\,12\,29  &  3\,13\,18  &  3\,13\,19  \\
     3\,15\,18  &  3\,15\,25  &  3\,17\,23  &  3\,20\,27  &  3\,21\,22  &  3\,22\,28  &  3\,23\,26  &  3\,26\,30  &  46\,13  &  47\,13  \\
     47\,22  &  48\,10  &  48\,20  &  4\,11\,14  &  4\,12\,14  &  4\,12\,17  &  4\,15\,24  &  4\,15\,27  &  4\,16\,25  &  4\,16\,26  \\
     4\,17\,21  &  4\,18\,19  &  4\,18\,27  &  4\,19\,25  &  4\,20\,28  &  4\,22\,29  &  4\,23\,24  &  4\,23\,30  &  4\,26\,28  &  4\,29\,30  \\
     56\,15  &  56\,17  &  579  &  57\,28  &  58\,15  &  59\,19  &  5\,10\,14  &  5\,10\,23  &  5\,11\,18  &  5\,11\,19  \\
     5\,12\,21  &  5\,12\,30  &  5\,13\,27  &  5\,14\,22  &  5\,16\,17  &  5\,18\,26  &  5\,20\,22  &  5\,20\,26  &  5\,21\,25  &  5\,23\,29  \\
     5\,24\,28  &  5\,25\,27  &  5\,29\,30  &  67\,14  &  67\,16  &  68\,10  &  68\,27  &  69\,28  &  69\,30  &  6\,10\,26  \\
     6\,11\,16  &  6\,12\,21  &  6\,12\,25  &  6\,13\,21  &  6\,14\,27  &  6\,15\,23  &  6\,17\,26  &  6\,18\,20  &  6\,18\,30  &  6\,19\,23  \\
     6\,19\,29  &  6\,22\,28  &  78\,16  &  78\,30  &  79\,11  &  7\,10\,15  &  7\,10\,25  &  7\,11\,26  &  7\,13\,20  &  7\,17\,20  \\
     7\,18\,21  &  7\,18\,23  &  7\,19\,24  &  7\,21\,30  &  7\,22\,27  &  7\,23\,25  &  7\,26\,27  &  7\,28\,29  &  89\,15  &  89\,17  \\
     8\,11\,17  &  8\,11\,22  &  8\,12\,23  &  8\,12\,25  &  8\,13\,18  &  8\,13\,21  &  8\,14\,29  &  8\,16\,21  &  8\,18\,25  &  8\,19\,24  \\
     8\,19\,26  &  8\,22\,26  &  8\,23\,24  &  8\,27\,28  &  8\,28\,30  &  9\,11\,29  &  9\,12\,18  &  9\,12\,26  &  9\,13\,20  &  9\,13\,23  \\
     9\,14\,22  &  9\,14\,25  &  9\,15\,29  &  9\,16\,27  &  9\,18\,22  &  9\,19\,21  &  9\,20\,30  &  9\,24\,26  &  9\,24\,27  &  9\,25\,28  \\
     10\,11\,21  &  10\,11\,27  &  10\,13\,28  &  10\,14\,24  &  10\,15\,16  &  10\,16\,19  &  10\,17\,20  &  10\,17\,29  &  10\,18\,21  &  10\,18\,25  \\
     10\,20\,22  &  10\,22\,30  &  10\,23\,27  &  10\,24\,30  &  10\,26\,29  &  11\,12\,13  &  11\,12\,16  &  11\,13\,28  &  11\,15\,20  &  11\,15\,24  \\
     11\,17\,19  &  11\,18\,20  &  11\,21\,24  &  11\,22\,29  &  11\,23\,26  &  11\,28\,30  &  12\,13\,23  &  12\,14\,24  &  12\,15\,20  &  12\,15\,28  \\
     12\,16\,29  &  12\,18\,19  &  12\,19\,27  &  12\,20\,30  &  12\,22\,26  &  12\,24\,27  &  13\,14\,26  &  13\,14\,30  &  13\,15\,29  &  13\,16\,22  \\
     13\,16\,24  &  13\,19\,27  &  13\,22\,25  &  13\,24\,26  &  13\,25\,29  &  14\,15\,21  &  14\,15\,30  &  14\,16\,18  &  14\,16\,28  &  14\,17\,18  \\
     14\,20\,21  &  14\,20\,25  &  14\,23\,28  &  14\,26\,29  &  15\,16\,26  &  15\,17\,19  &  15\,17\,27  &  15\,18\,26  &  15\,19\,30  &  15\,21\,23  \\
     15\,25\,28  &  16\,17\,21  &  16\,18\,30  &  16\,19\,20  &  16\,20\,24  &  16\,25\,30  &  16\,28\,29  &  17\,18\,28  &  17\,22\,23  &  17\,22\,27  \\
     17\,24\,29  &  17\,24\,30  &  17\,25\,30  &  17\,26\,28  &  18\,22\,24  &  18\,23\,29  &  18\,27\,28  &  19\,20\,21  &  19\,22\,23  &  19\,22\,30  \\
     19\,25\,26  &  20\,23\,25  &  20\,23\,28  &  20\,26\,27  &  21\,23\,30  &  21\,24\,29  &  21\,25\,26  &  21\,27\,29  &  22\,24\,25  &  25\,27\,29.
\end{tabular*}
}

\vspace{3mm}

\mbox{}

\section{A generalization of Ringel's cyclic \mathversion{bold}$7{\rm mod}12$ series\mathversion{normal}}
\label{sec:generalized_Ringel}

Ringel's cyclic $7{\rm mod}12$ series \cite{Ringel1961} provides for every $k\geq 0$ a neighborly triangulation 
with \mbox{$n=7+12k$} vertices of the orientable surface of genus $g=12k^2+7k+1$. We generalize Ringel's series
to a two-parameter series of cyclic triangulations.

\begin{thm}\label{thm_generalized_Ringel}
For every $k\geq 0$ and every $n\geq 7+12k$ the orientable surface of genus $g=kn+1$ has a triangulation $R_{k,n}$ 
with $n$ vertices and vertex-transitive cyclic symmetry. The examples $R_{k,n}$ are equivelar of type $\{3,6+12k\}$ 
and, in terms of $q=6+12k$, the examples have $f$-vector $f=(n,3(1+2k)n,2(1+2k)n)=(n,\frac{nq}{2},\frac{nq}{3})$
and genus $g=1+\frac{(q-6)n}{12}$.
\end{thm}

\textbf{Proof:} Step I: Let $k\geq 0$ and $n=7+12k$. Then $f=(n,3(1+2k)n,2(1+2k)n)=(n,\binom{n}{2},\frac{2}{3}\binom{n}{2})$,
which is the $f$-vector of a neighborly triangulation of a (orientable or non-orientable) surface 
of Euler characteristic $\chi=n-\frac{1}{3}\binom{n}{2}$. Ringel's cyclic $7{\rm mod}12$ series \cite{Ringel1961}, \cite[Ch.~2.3]{Ringel1974}
provides examples $R_{k,7+12k}$ of triangulated orientable surfaces with the desired parameters. 

\begin{figure}
\begin{center}
\scriptsize
\psfrag{1}{1}
\psfrag{2}{2}
\psfrag{3}{3}
\psfrag{4}{4}
\psfrag{5}{5}
\psfrag{6}{6}
\psfrag{2k-5}{$2k-5$}
\psfrag{2k-4}{$2k-4$}
\psfrag{2k-3}{$2k-3$}
\psfrag{2k-2}{$2k-2$}
\psfrag{2k-1}{$2k-1$}
\psfrag{2k}{$2k$}
\psfrag{2k+1}{$2k+1$}
\psfrag{2k+2}{$2k+2$}
\psfrag{2k+3}{$2k+3$}
\psfrag{2k+4}{$2k+4$}
\psfrag{3k}{$3k$}
\psfrag{3k+1}{$3k+1$}
\psfrag{3k+2}{$3k+2$}
\psfrag{3k+3}{$3k+3$}
\psfrag{3k+4}{$3k+4$}
\psfrag{3k+5}{$3k+5$}
\psfrag{4k}{$4k$}
\psfrag{4k+1}{$4k+1$}
\psfrag{4k+2}{$4k+2$}
\psfrag{4k+3}{$4k+3$}
\psfrag{4k+4}{$4k+4$}
\psfrag{4k+5}{$4k+5$}
\psfrag{5k+1}{$5k+1$}
\psfrag{5k+2}{$5k+2$}
\psfrag{5k+3}{$5k+3$}
\psfrag{5k+4}{$5k+4$}
\psfrag{5k+5}{$5k+5$}
\psfrag{5k+6}{$5k+6$}
\psfrag{6k+1}{$6k+1$}
\psfrag{6k+2}{$6k+2$}
\psfrag{6k+3}{$6k+3$}
\includegraphics[width=\linewidth]{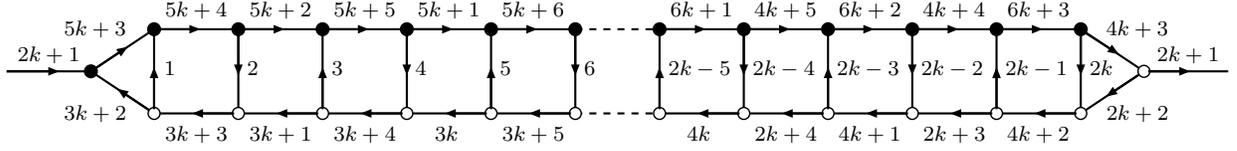}
\end{center}
\caption{The directed graph $G_{k,7+12k}$.}
\label{fig:7mod12}
\end{figure}

We follow Ringel's description \cite[pp.~25--28]{Ringel1974} of the series $R_{k,7+12k}$ and consider a directed graph $G_{k,7+12k}$ 
for any given $k\geq 0$ as displayed in Figure~\ref{fig:7mod12}, where the leftmost arc and the rightmost arc are to be identified. 
The $6k+3$ arcs of the digraph are labeled by $1,\dots,6k+3$ in a way, such that Kirchhoff's rule holds at every vertex, 
i.e., at each vertex the sum of the values of the ingoing arcs equals the sum of the values of the outgoing arcs. 
The elements $1,\dots,6k+3$ are regarded as elements of the cyclic group ${\mathbb Z}_n={\mathbb Z}_{7+12k}$ 
with inverses $\overline{1},\dots,\overline{6k+3}$ and neutral element $0$. We let the vertices of $R_{k,7+12k}$ 
be the elements of ${\mathbb Z}_{7+12k}$. Next, we read off the vertex-link of the vertex $0$ in the triangulation $R_{k,7+12k}$ 
from Figure~\ref{fig:7mod12} by ``walking along'' the arcs of the digraph:  Start with the arc with label, say, $1$ and follow the arc 
in the indicated direction. At a vertex, turn left whenever the vertex is black and turn right whenever the vertex is white, respectively. 
Thus, at the end of arc~$1$ we turn left and follow the arc with label $5k+3$ in the opposite direction of the arc. Then we turn left
again and follow the arc $3k+2$ in opposite direction, etc. Along our walk we record the labels of the arcs, respectively their 
inverses, if an arc is traversed in opposite direction. Thus, we record (for the link of the vertex $0$):
\begin{center}
\scriptsize
\begin{tabular}{c@{\hspace{3mm}}l@{\hspace{1mm}}l@{\hspace{1mm}}l@{\hspace{1mm}}l@{\hspace{1mm}}l@{\hspace{1mm}}l@{\hspace{1mm}}l@{\hspace{1mm}}l@{\hspace{1mm}}l@{\hspace{1mm}}l@{\hspace{1mm}}l@{\hspace{1mm}}l}
  0: & 1, & $\overline{5k+3}$, & $\overline{3k+2}$, & $\overline{3k+3}$, & \dots, & $\overline{2k+2}$, & $2k+1$, & $5k+3$, & $5k+4$, & \dots & 2, & $3k+3$.
\end{tabular}
\end{center}
Every arc of the digraph is traversed twice, once into the indicated direction and once into the opposite direction.
Hence, the link of the vertex $0$ is a circle with $2\cdot(6k+3)=6+12k$ vertices. 
For convenience, we write, from now on, $n-a$ for the inverse $\overline{a}$ of an element $a\in{\mathbb Z}_n$. 

\pagebreak

The vertex links of the vertices $1$, \dots, $n-1$ are obtained from the link of $0$ by cyclic shifts:

\begin{center}
\scriptsize
\begin{tabular}{c@{\hspace{3mm}}l@{\hspace{1mm}}l@{\hspace{1mm}}l@{\hspace{1mm}}l@{\hspace{1mm}}l@{\hspace{1mm}}l@{\hspace{1mm}}l@{\hspace{1mm}}l@{\hspace{1mm}}l@{\hspace{1mm}}l@{\hspace{1mm}}l@{\hspace{1mm}}l}
  0: &  1, & $n-5k-3$, & $n-3k-2$, & $n-3k-3$, & \dots, & $n-2k-2$, & $2k+1$, & $5k+3$, & $5k+4$, & \dots, & 2, & $3k+3$\\
  1: &  2, & $n-5k-2$, & $n-3k-1$, & $n-3k-2$, & \dots, & $n-2k-1$, & $2k+2$, & $5k+4$, & $5k+5$, & \dots, & 3, & $3k+4$\\
  2: &  3, & $n-5k-1$, & $n-3k$,   & $n-3k-1$, & \dots, & $n-2k$,   & $2k+3$, & $5k+5$, & $5k+6$, & \dots, & 4, & $3k+5$\\
\vdots & & & & & \vdots & & & & & \vdots & & 
\end{tabular}
\end{center}

Step II: In order to describe a triangulated surface, we can, instead of the vertex links, give the triangles of the surface.
In the case of the Ringel series $R_{k,7+12k}$ with cyclic symmetry it is sufficient to give one generating triangle
from each of the $2+4k$ orbits of triangles, with each orbit of size $7+12k$. We can even choose (and this will be important
in Step III) a set of generating triangles with vertices from the subset $\{0,1,2,\dots,6k+3\}$; see Table~\ref{tbl:generators}.

\begin{table}
\begin{center}
\scriptsize
\begin{tabular}{c@{\hspace{7mm}}c}
$[0,1,3k+3]$, & $[0,2,3k+3]$, \\ 
$[0,3,3k+4]$, & $[0,4,3k+4]$, \\ 
$[0,5,3k+5]$, & $[0,6,3k+5]$, \\ 
\vdots & \vdots \\
$[0,2k-3,4k+1]$, & $[0,2k-2,4k+1]$, \\ 
$[0,2k-1,4k+2]$, & $[0,2k,4k+2]$, \\[4mm] 
$[0,4k+3,6k+3]$, & $[0,2k+2,4k+3]$, \\[4mm] 
$[0,4k+4,6k+2]$, & $[0,4k+4,6k+3]$, \\ 
$[0,4k+5,6k+1]$, & $[0,4k+5,6k+2]$, \\ 
\vdots & \vdots \\
$[0,5k+1,5k+5]$, & $[0,5k+1,5k+6]$, \\ 
$[0,5k+2,5k+4]$, & $[0,5k+2,5k+5]$, \\[4mm] 
$[0,2k+1,5k+3]$, & $[0,5k+3,5k+4]$.
\end{tabular}
\end{center}
\caption{The set of generating triangles for $R_{k,n}$.}
\label{tbl:generators}
\end{table}

In each orbit there are three triangles that contain the vertex $0$ and thus contribute an edge to the vertex-link
of $0$. For the choice of generating triangles in Table~\ref{tbl:generators}, these $3\cdot(2+4k)$ edges are as follows:

\begin{center}
\scriptsize
\begin{tabular}{@{}c@{\hspace{1mm}}c@{\hspace{1mm}}c@{\hspace{4mm}}c@{\hspace{1mm}}c@{\hspace{1mm}}c@{}}
$[1,3k+3]$, & $[3k+2,n-1]$, & $[n-3k-3,n-3k-2]$, & $[2,3k+3]$, & $[3k+1,n-2]$, & $[n-3k-3,n-3k-1]$, \\ 
$[3,3k+4]$, & $[3k+1,n-3]$, & $[n-3k-4,n-3k-1]$, & $[4,3k+4]$, & $[3k,n-4]$,   & $[n-3k-4,n-3k]$, \\ 
$[5,3k+5]$, & $[3k,n-5]$,   & $[n-3k-5,n-3k]$,   & $[6,3k+5]$, & $[3k-1,n-6]$, & $[n-3k-5,n-3k+1]$, \\ 
\vdots & \vdots & \vdots & \vdots & \vdots & \vdots \\
$[2k-3,4k+1]$, & $[2k+4,n-2k+3]$, & $[n-4k-1,n-2k-4]$, & $[2k-2,4k+1]$, & $[2k+3,n-2k+2]$, & $[n-4k-1,n-2k-3]$, \\ 
$[2k-1,4k+2]$, & $[2k+3,n-2k+1]$, & $[n-4k-2,n-2k-3]$, & $[2k,4k+2]$,   & $[2k+2,n-2k]$,   & $[n-4k-2,n-2k-2]$, \\[4mm] 
$[4k+3,6k+3]$, & $[2k,n-4k-3]$,   & $[n-6k-3,n-2k]$,   & $[2k+2,4k+3]$, & $[2k+1,n-2k-2]$, & $[n-4k-3,n-2k-1]$, \\[4mm] 
$[4k+4,6k+2]$, & $[2k-2,n-4k-4]$, & $[n-6k-2,n-2k+2]$, & $[4k+4,6k+3]$, & $[2k-1,n-4k-4]$, & $[n-6k-3,n-2k+1]$, \\ 
$[4k+5,6k+1]$, & $[2k-4,n-4k-5]$, & $[n-6k-1,n-2k+4]$, & $[4k+5,6k+2]$, & $[2k-3,n-4k-5]$, & $[n-6k-2,n-2k+3]$, \\ 
\vdots & \vdots & \vdots & \vdots & \vdots & \vdots \\
$[5k+1,5k+5]$, & $[4,n-5k-1]$, & $[n-5k-5,n-4]$, & $[5k+1,5k+6]$, & $[5,n-5k-1]$, & $[n-5k-6,n-5]$, \\ 
$[5k+2,5k+4]$, & $[2,n-5k-2]$, & $[n-5k-4,n-2]$, & $[5k+2,5k+5]$, & $[3,n-5k-2]$, & $[n-5k-5,n-3]$, \\[4mm] 
$[2k+1,5k+3]$, & $[3k+2,n-2k-1]$, & $[n-5k-3,n-3k-2]$, & $[5k+3,5k+4]$, & $[1,n-5k-3]$, & $[n-5k-4,n-1]$.
\end{tabular}
\end{center}

Of course, both link-descriptions yield the same triangulated circle as the link of the vertex~$0$. 
In Ringel's description, the edges of the link are read off from the digraph $G_{k,7+12k}$ in consecutive order, 
whereas the description via the set of generating triangles lists the edges of the link in correspondence 
to the orbits of triangles in which they occur.

\medskip

\emph{Remark:} The labeled graph $G_{k,7+12k}$ degenerates for $k=0$ and $k=1$, and so do the orbit descriptions of $R_{\,0,n}$ and $R_{\,1,n}$. 
For $k=1$, we have the orbit generating triangles
\begin{center}
\small
\begin{tabular}{c@{\hspace{7mm}}c}
$[0,1,6]$, & $[0,2,6]$, \\ 
$[0,7,9]$, & $[0,4,7]$, \\ 
 --- &  --- \\
$[0,3,8]$, & $[0,8,9]$,
\end{tabular}
\end{center}
whereas for $k=0$, the orbit generating triangles are
\begin{center}
\small
\begin{tabular}{c@{\hspace{7mm}}c}
$[0,1,3]$, & $[0,2,3]$,
\end{tabular}
\end{center}
which precisely are the orbit generating triangles of Altshuler's series of cyclic tori $T^2(n)$ for $n\geq 7$;
see Figure~\ref{fig:T_2_n}.

\medskip

\begin{figure}
\begin{center}
\psfrag{i}{$i$}
\psfrag{j}{$j$}
\psfrag{k}{$k$}
\psfrag{l}{$l$}
\includegraphics[width=.425\linewidth]{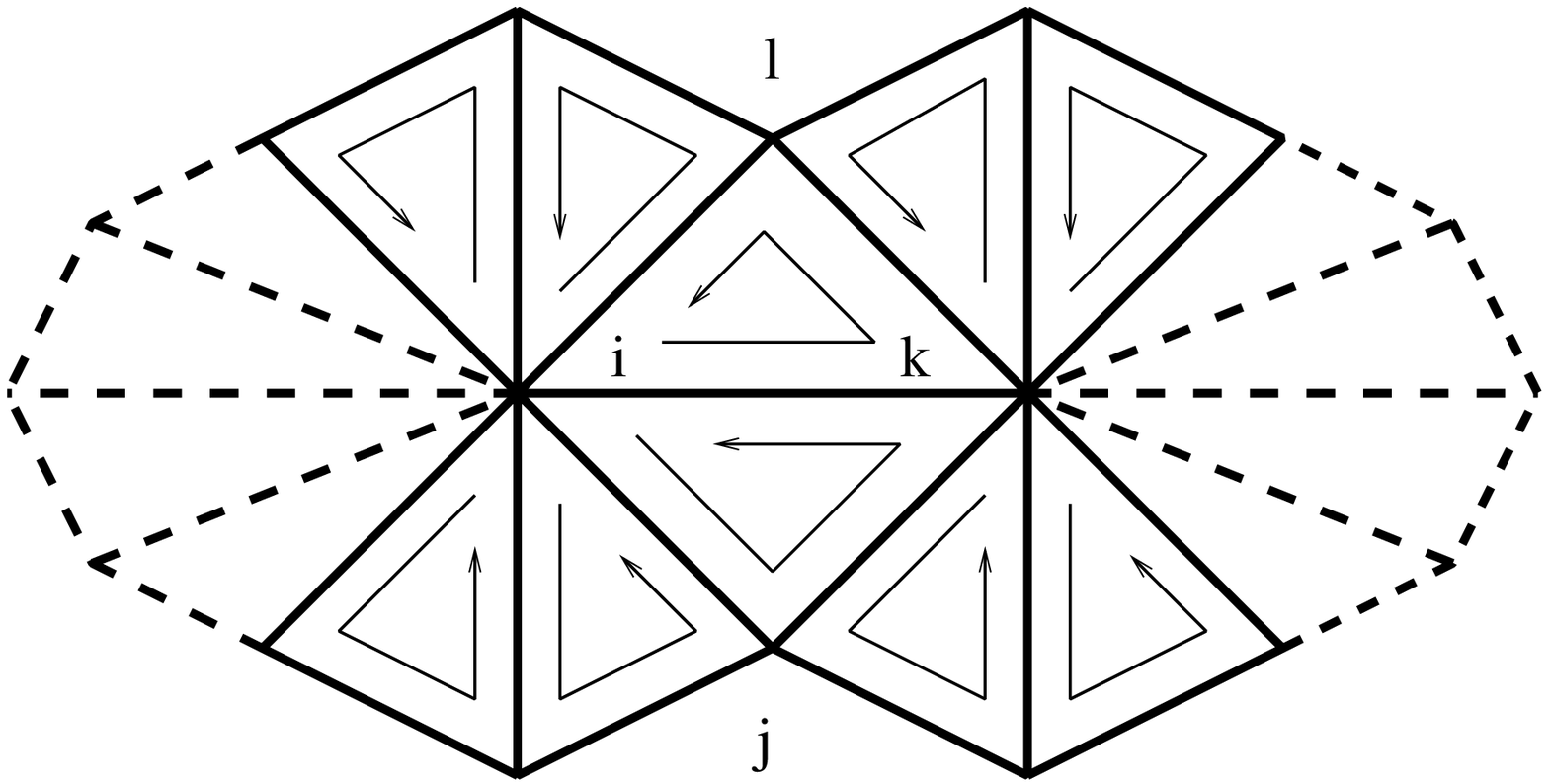}
\end{center}
\caption{Rule $R^*$.}
\label{fig:rule_R_star}
\end{figure}

Step III: For $n\geq 7+12k$, let the vertices of $R_{k,n}$ be $\{0,1,2,\dots,n-1\}$ and define $R_{k,n}$
via the same set of $2+4k$ generating triangles from Table~\ref{tbl:generators} 
with respect to the action of the cyclic group~${\mathbb Z}_n$.
Since all vertices of the generating triangles are elements of the set $\{0,1,2,\dots,6k+3\}$,
the link of vertex $0$ has the same description as before in the case $n=7+12k$.
In other words, the `symbol' $n=7+12k$ in the link of $0$ in $R_{k,7+12k}$ is replaced by the `symbol' $n\geq7+12k$ in the link of $0$ in $R_{k,n}$.
This way, no unwanted identifications of vertices in the link take place, so the link is one circle of length $6+12k$.
It follows, that $R_{k,n}$ is a triangulated surface, which is equivelar of type $\{3,6+12k\}$, has $f$-vector $f=(n,3(1+2k)n,2(1+2k)n)$,
and genus $g=kn+1$.

It remains to show that $R_{k,n}$ is orientable for all $n\geq 7+12k$, $k\geq 0$. In Ringel's description of the link of $0$
in the case $n=7+12k$, $k\geq 0$, the edges of the link are listed in consecutive order, which induces an orientation 
of the triangles of the star of the vertex $0$. Since the stars of the other vertices are obtained from the star of $0$ by cyclic shifts,
all vertex-stars inherit an orientation from the orientation of the star of $0$. For $n=7+12k$, the orientations of the
individual vertex-stars are compatible, as is expressed by Ringel's rule $R^*$ for the vertex-links \cite[p.~28]{Ringel1974}:

\medskip

\emph{Rule $R^*$}:\, If in row $i$ one has\, $i:\,\dots jkl \dots$, then row $k$ appears as\, $k:\,\dots lij \dots $ (see Figure~\ref{fig:rule_R_star}). 

\medskip

In particular, rule $R^*$ holds for row $0$, that is, if we have\, $0:\,\dots jkl \dots$, 
then row $k$ appears as\, $k:\,\dots l0j \dots $. If we apply a cyclic shift to row $k$ by adding $-k$,
we obtain\, $0:\,\dots l-k,-k,j-k \dots $. 

For $R_{k,7+12k}$, the link of $0$ has the set $P\cup N$ as its set of vertices, with $P=\{1,2,\dots,6k+3\}$ and $N=\{n-1,n-2,\dots,n-(6k+3)\}$.
In the transition from $R_{k,7+12k}$ to $R_{k,n}$, that is, from $n=7+12k$ to $n\geq 7+12k$, 
a vertex $a\in P$ of the link of $0$ remains unchanged, whereas a vertex $b$ in the link of $0$ in $R_{k,7+12k}$ with
$b\in N$ can be written as $b=(7+12k)-r$ with $r\in P$, and is replaced by $n-r$.
As a consequence, if we have\, $0:\,\dots jkl \dots$ for the row $0$ in $R_{k,n}$ for $n\geq 7+12k$, 
then also\, $0:\,\dots l-k,-k,j-k \dots $, as originally in $R_{k,7+12k}$.
If we now add $k$ to this row, we obtain\, $k:\,\dots l0j \dots $, so the orientation of the star of $0$
is compatible with the orientations of all the stars of the vertices in the link of $0$.
By the cyclic symmetry it follows that the orientations of any two neighboring vertex-stars are compatible.
Therefore, $R_{k,n}$ is orientable.\hfill$\Box$

\medskip
\smallskip

The generalized Ringel series $R_{k,n}$ interpolates between Ringel's cyclic neighborly series $R_{k,7+12k}$ 
and Altshuler's cyclic series of $6$-equivelar tori $T^2(n)=R_{0,n}$. The examples $R_{k,n}$ settle 
Conjecture~\ref{conj:triples} for $q=6+12k$, $k\geq 1$, for the orientable surfaces with $\chi=-2kn$.

\begin{cor}
The two subseries $R_{k,7+12k+1}$ and $R_{k,7+12k+2}$ of the generalized Ringel series  
provide for $k\geq 1$ infinitely many examples of non-neighborly cyclic triangulations
with $q=6+12k=\Bigl\lfloor\tfrac{1}{2}(5+\sqrt{49-24\chi})\Bigl\rfloor$,
i.e., with equality in~(\ref{eq:bound_for_q}). By subdivision, these examples yield
$2q$-covered triangulations with equality in~(\ref{eq:bound_for_d}).
\end{cor}

\textbf{Proof:} Let for $k\geq 1$, $n=7+12k+m$ with $m\geq 0$. By Theorem~\ref{thm_generalized_Ringel},
the examples $R_{k,n}$ of the generalized Ringel series have Euler characteristic $\chi=-2kn$.

For $m=0$, $k\geq 1$, we have $n=7+12k$. Therefore, the example $R_{k,7+12k}$ is neighborly
with $q=\Bigl\lfloor\tfrac{1}{2}(5+\sqrt{49-24\chi})\Bigl\rfloor=\Bigl\lfloor\tfrac{1}{2}(5+\sqrt{49+24\cdot 2k(7+12k)})\Bigl\rfloor=\Bigl\lfloor\tfrac{1}{2}(5+\sqrt{(7+24k)^2})\Bigl\rfloor=6+12k$.

We next determine those $m\geq 1$, such that $\Bigl\lfloor\tfrac{1}{2}(5+\sqrt{49+24\cdot 2k(7+12k+m)})\Bigl\rfloor=6+12k$,
which is equivalent to find those $m\geq 1$, for which $\sqrt{(7+24k)^2+48km}<(7+24k)+2$,
that is, $m<2+\frac{2}{3k}$. Thus, $1\leq m\leq 2$.\hfill$\Box$

\begin{cor}
For $k\geq 1$, the series $R_{k,7+12k+1}$ provides vertex-minimal triangulations of the orientable
surfaces of Euler characteristic $\chi=-2kn$ with vertex-transitive cyclic symmetry on $n=7+12k+1$ vertices.
\end{cor}

\textbf{Proof:} For $k\geq 1$, $\Bigl\lceil\tfrac{1}{2}(7+\sqrt{49+24\chi})\Bigl\rceil=\Bigl\lceil\tfrac{1}{2}(7+\sqrt{49+24\cdot 2k(7+12k+1)})\Bigl\rceil=7+12k+1$.\linebreak
\mbox{}\hfill$\Box$

\bibliography{.}

\bigskip
\bigskip
\smallskip

\small

\noindent
Frank H. Lutz\\
Institut f\"ur Mathematik\\
Technische Universit\"at Berlin\\
Stra\ss e des 17.\ Juni 136\\
10623 Berlin, Germany\\
{\tt lutz@math.tu-berlin.de}

\vspace{6mm}

\noindent
Thom Sulanke\\
Department of Physics\\
Indiana University\\
Bloomington, Indiana 47405, USA\\
{\tt tsulanke@indiana.edu}

\vspace{6mm}

\noindent
Ashish K.~Upadhyay,\, Anand K.~Tiwari\\
Department of Mathematics\\
Indian Institute of Technology Patna\\ 
Patliputra Colony, Patna -- 800\,013, India\\
{\tt \{upadhyay,anand\}@iitp.ac.in}

\end{document}